\documentclass{amsart}
\usepackage[all]{xy}

\newtheorem{theorem}{Theorem}[section]
\newtheorem{definition}[theorem]{Definition}
\newtheorem{proposition}[theorem]{Proposition}

\newtheorem{lemma}[theorem]{Lemma}
\newtheorem{question}[theorem]{Question}

\theoremstyle{remark}
\newtheorem*{remark}{Remark}

\newcommand{\GF}{\mathbb{F}}
\newcommand{\BH}{\mathbb{H}}

\newcommand{\BP}{\mathbb{P}}
\newcommand{\Qbar}{\overline{\mathbb{Q}}}
\newcommand{\Bmu}{\mathbb{\mu}}
\newcommand{\Q}{\mathbb{Q}}

\newcommand{\Z}{\mathbb{Z}}

\newcommand{\calA}{\mathcal{A}}

\newcommand{\calL}{\mathcal{L}}
\newcommand{\calO}{\mathcal{O}}

\newcommand{\frakA}{\mathfrak{A}}
\newcommand{\frakB}{\mathfrak{B}}
\newcommand{\frakp}{\mathfrak{p}}
\newcommand{\frakP}{\mathfrak{P}}

\newcommand{\mybar}[1]{#1\llap{$\overline{\phantom{\rm#1}}$}}

\newcommand{\GFbar}{\mybar{\mathbb{F}}}

\newcommand{\kbar}{\mybar{k}}
\newcommand{\xbar}{\mybar{x}}

\newcommand{\omegabar}{\mybar{\omega}}

\newcommand{\eps}{\varepsilon}
\renewcommand{\epsilon}{\varepsilon}
\renewcommand{\phi}{\varphi}
\renewcommand{\theta}{\vartheta}

\newcommand\as[1]{\renewcommand\arraystretch{#1}}
\newcommand\et{E\times E}

\newcommand\ff[1]{\GF_{q^{#1}}}
\newcommand\fq{\GF_q}
\newcommand\fp{\GF_p}
\newcommand\fpt{\GF_{p^2}}

\newcommand\gl[2]{\GL_{#1}(#2)}
\renewcommand{\hat}{\widehat}
\newcommand\nprc{\text{\rm nprinc}}
\newcommand\pgl[2]{\PGL_{#1}(#2)}
\newcommand\pr[1]{\BP^{#1}}
\newcommand\prc{\text{\rm princ}}
\newcommand\sii{\,\Longleftrightarrow\,}
\newcommand\tq{\,\,|\,\,}

\DeclareMathOperator{\Aut}{Aut}

\DeclareMathOperator{\End}{End}
\DeclareMathOperator{\Gal}{Gal}
\DeclareMathOperator{\GL}{GL}
\DeclareMathOperator{\gra}{Graph}
\DeclareMathOperator{\Hom}{Hom}

\DeclareMathOperator{\Jac}{Jac}
\DeclareMathOperator{\Norm}{Norm}
\DeclareMathOperator{\NS}{NS}

\DeclareMathOperator{\PGL}{PGL}

\DeclareMathOperator{\SL}{SL}

\newcommand\lra{\longrightarrow}
\newcommand\llra{\relbar\joinrel\longrightarrow}

\newcommand{\mapright}[1]{\mathop{\longrightarrow}\limits^{#1}}

\newcommand\longmapright[1]{\mathop{\llra}\limits^{#1}}

\begin{document}

\title[Jacobians in isogeny classes]
{Jacobians in isogeny classes of abelian surfaces over finite fields}

\author{Everett W.~Howe}
\address{Center for Communications Research,
         4320 Westerra Court,
         San Diego, CA 92121-1967, USA.}
\email{however@alumni.caltech.edu}
\urladdr{http://www.alumni.caltech.edu/~however/}

\author{Enric Nart}
\address{Departament de Matem\`atiques,
        Universitat Aut\`onoma de Barcelona,
        Edifici C, 08193 Bellaterra, Barcelona, Spain.}
\email{nart@mat.uab.es}

\author{Christophe Ritzenthaler}
\address{Institut de Math\'ematiques de Luminy,
         UMR 6206 du CNRS,
         Luminy, Case 907, 13288 Marseille, France.}
\email{ritzenth@iml.univ-mrs.fr}

\thanks{The second author was supported in
part by project MTM2006-11391 from the Spanish MEC}

\date{2 April 2007}



\begin{abstract}
We give a complete answer to the question of which polynomials occur
as the characteristic polynomials of Frobenius for genus-$2$ curves
over finite fields.
\end{abstract}

\maketitle

\section{Introduction}
\label{S:paperintro}
The \emph{Weil polynomial} of an abelian variety over a finite field $\fq$
is the characteristic polynomial of its Frobenius endomorphism; the
\emph{Weil polynomial} of a curve over $\fq$ is the Weil polynomial of its
Jacobian.  In this paper we determine which polynomials occur as the Weil
polynomials of genus-$2$ curves over a finite field.

Weil's `Riemann Hypothesis' shows that the Weil polynomial of an abelian
surface over $\fq$ has the form
$$x^4 + ax^3 + bx^2 + aqx + q^2,$$
and the Honda-Tate theorem~\cite{tate}  makes it a straightforward matter to
determine which such polynomials come from abelian surfaces
(see \cite[Thm.~1.1]{ruc}, \cite[Thm.~2.9]{mn}, and the Appendix to this
paper).  Since two abelian varieties over $\fq$ are isogenous to one another
if and only if they have the same Weil polynomial~\cite{tate-invent},
we may phrase our main question as follows.

\begin{question}
\label{Q:main}
Suppose $f = x^4+ax^3+bx^2+aqx+q^2$ is the Weil polynomial for an isogeny
class of abelian surfaces over $\fq$. Is there a projective smooth genus-$2$
curve over $\fq$ whose Weil polynomial is equal to $f$?
\end{question}

While this question has been settled in many special cases (as we
explain in detail below), until now there have been two kinds of isogeny
classes for which the question has remained largely unanswered:
the split isogeny classes and the supersingular isogeny classes.
We analyze these two remaining cases and provide
a complete answer to Question~\ref{Q:main}.

\begin{theorem}
\label{T:main}
Let $f = x^4 +ax^3+bx^2+aqx+q^2$ be the Weil polynomial of an isogeny
class $\calA$ of abelian surfaces over $\fq$, where $q$ is a power
of a prime $p$.
\begin{enumerate}
\item
Suppose that $\calA$ contains a product of elliptic curves,
so that $f$ can be written as a product
$$f = (x^2 - sx + q) (x^2 - tx + q)$$
where the two factors are the Weil polynomials of isogeny classes of
elliptic curves over $\fq$ and where we may assume that $|s|\ge |t|$.
Then $\calA$ does not contain a Jacobian if and only if the conditions
in one of the rows of Table~\ref{Table:MainSplit} are met.
\item
Suppose that $\calA$ is simple.  Then $\calA$ does not contain a
Jacobian if and only if the conditions
in one of the rows of Table~\ref{Table:MainSimple} are met.
\end{enumerate}
\end{theorem}

\begin{table}
\begin{center}
\as{1.1}
\begin{tabular}{|c|l|l|}
\hline
$p$-rank of $\calA$ & Condition on $p$ and $q$ & Conditions on $s$ and $t$\\
\hline\hline
 --- & ---                     & $|s-t| = 1$                           \\ \hline
 $2$ & ---                     & $s=t$ and $t^2 - 4q\in \{-3,-4,-7\}$  \\ \cline{2-3}
     & $q=2$                   & $|s|=|t|=1$ and $s\neq t$             \\ \hline
 $1$ & $q$ square              & $s^2 = 4q$ and $s-t$ squarefree       \\ \hline
     & $p>3$                   & $s^2 \neq t^2$                        \\ \cline{2-3}
     & $p=3$ and $q$ nonsquare & $s^2 = t^2 = 3q$                      \\ \cline{2-3}
 $0$ & $p=3$ and $q$ square    & $s-t$ is not divisible by $3\sqrt{q}$ \\ \cline{2-3}
     & $p=2$                   & $s^2 - t^2$ is not divisible by $2q$  \\ \cline{2-3}
     & $q=2$ or $q=3$          & $s = t$                               \\ \cline{2-3}
     & $q=4$ or $q=9$          & $s^2 = t^2 = 4q$                      \\ \hline
\end{tabular}
\end{center}
\vspace{1ex}
\caption{Conditions that ensure that the split isogeny class
with Weil polynomial $(x^2 - sx + q)(x^2 - tx + q)$
does not contain a Jacobian.   Here we assume that $|s|\ge |t|$.}
\label{Table:MainSplit}
\end{table}

\begin{table}
\begin{center}
\as{1.1}
\begin{tabular}{|c|l|l|}
\hline
$p$-rank of $\calA$ & Condition on $p$ and $q$ & Conditions on $a$ and $b$\\
\hline\hline
--- & ---                                 & $a^2-b = q$ and $b<0$ and \\
    &                                     & all prime divisors of $b$ are $1\bmod3$\\ \hline
$2$ & ---                                 & $a = 0$ and $b = 1-2q$ \\ \cline{2-3}
    & $p>2$                               & $a = 0$ and $b = 2-2q$ \\ \hline
    & $p\equiv 11\bmod 12$ and $q$ square & $a = 0$ and $b = -q$   \\ \cline{2-3}
$0$ & $p= 3$ and $q$ square               & $a = 0$ and $b = -q$   \\ \cline{2-3}
    & $p= 2$ and $q$ nonsquare            & $a = 0$ and $b = -q$   \\ \cline{2-3}
    & $q=2$ or $q=3$                      & $a = 0$ and $b = -2q$  \\ \hline
\end{tabular}
\end{center}
\vspace{1ex}
\caption{Conditions that ensure that the simple isogeny class with
Weil polynomial $x^4 + ax^3 + bx^2 + aqx + q^2$
does not contain a Jacobian.}
\label{Table:MainSimple}
\end{table}

The analog of Question~\ref{Q:main} for elliptic curves was answered by
Waterhouse~\cite[Thm.~4.1]{wa}.  The generalization from elliptic
curves to genus-$2$ curves is surely quite natural, but to the best of our
knowledge Question~\ref{Q:main} did not occur in print until 1990,
when R\"uck~\cite{ruc} provided some sufficient conditions for a positive
answer to the question.  In the literature starting with R\"uck we find a large
variety of methods and techniques that provide both positive and negative
answers to Question~\ref{Q:main} for particular classes of Weil polynomials.
Almost all of the positive results are based on the following theorem of
Weil; the version we give here is due to Gonz\'alez, Gu\`ardia, and
Rotger~\cite[Thm.~3.1]{GGR}.

\begin{theorem}[Weil]
\label{T:Weil}
Let $(A,\lambda)$ be a principally polarized abelian surface defined over
a field $k$.  Then $(A,\lambda)$ is either
\begin{itemize}
\item[{\rm (a)}] the polarized Jacobian of a genus-$2$ curve over $k$,
\item[{\rm (b)}] the product of two polarized elliptic curves over $k$, or
\item[{\rm (c)}] the restriction of scalars of a polarized elliptic curve
                 over a quadratic extension of~$k$.
\end{itemize}
Furthermore, these three possibilities are mutually exclusive.
\end{theorem}

Let us say that an isogeny class of abelian varieties is
\emph{principally polarizable} if it contains a principally polarized variety.
In light of Weil's theorem, if an isogeny class of abelian surfaces over $\fq$
is simple over~$\ff2$, then it contains a Jacobian if and only if it is
principally polarizable.  The problem of determining the principally
polarizable isogeny classes of abelian varieties was studied by the first
author in a series of papers \cite{howe1}, \cite{howe2}, \cite{howe3}, where he
expressed the obstruction to the existence of principal polarizations in terms
of the vanishing of an element of a group constructed from the Grothendieck
group of the category of finite group schemes that can be embedded in
varieties in the isogeny class.
Recall that an abelian surface over a field of characteristic $p>0$
is said to be \emph{ordinary} when its $p$-rank is $2$,
\emph{mixed} when its $p$-rank is $1$, and
\emph{supersingular} when its $p$-rank is $0$;
the $p$-rank of an abelian surface over a finite field
can be read from the Newton polygon of the Weil polynomial in a
well-known way.
The principally polarizable isogeny classes of ordinary abelian
surfaces over a finite field were determined
in~\cite[Thm.~1.3]{howe1}. In \cite{mn} it was proved that all mixed isogeny
classes over a finite field are principally polarizable by applying criteria
developed in \cite{howe2}; in particular, since the simple mixed classes are
absolutely simple they always contain Jacobians. Finally, the supersingular
case was worked out in \cite{hmnr} using the ideas of \cite{howe2}
and~\cite{howe3}.
Gathering all these results one obtains:

\begin{theorem}[\cite{hmnr}]
\label{T:PPAS}
Let $\calA$ be an isogeny class of abelian surfaces over $\fq$ with Weil
polynomial $x^4 + ax^3 + bx^2 + aqx + q^2$. Then $\calA$ is not principally
polarizable if and only if the following three conditions are satisfied{\/\rm:}
\begin{itemize}
\item[{\rm (a)}]  $a^2-b=q$,
\item[{\rm (b)}]  $b<0$, and
\item[{\rm (c)}]  all prime divisors of $b$ are congruent to $1$ modulo $3$.
\end{itemize}
\end{theorem}

This result, together with Weil's theorem, answers Question~\ref{Q:main} for
every isogeny class that is simple over~$\ff2$.

The answer to Question~\ref{Q:main} for the simple ordinary isogeny classes
that split over $\ff2$ was determined by the first author and Maisner.
The Weil polynomial of such an isogeny class is of the form $x^4 + bx^2 + q^2$.
The first author~\cite[App.]{mn} proved that when $b = 1 - 2q$ there is no
curve with the given Weil polynomial by showing that such a curve would have an
automorphism whose existence is incompatible with the number of rational
points on the curve over~$\ff2$.  For $p>2$, the first author~\cite{howecomp}
used a counting argument to show that when $b = 2-2q$ there is again no curve
with the given Weil polynomial.  He found explicit formulas for the
number of principally polarized surfaces $(A,\lambda)$ with $A$ belonging to
the given isogeny class, as well as for the number of these polarized surfaces
that are restrictions of scalars of elliptic curves over~$\ff2$.  The formulas
for these two numbers both involve arithmetic invariants of the biquadratic
field generated by the Weil polynomial, and a comparison of the
two numbers using the Brauer relations shows that they coincide;
thus, all $(A,\lambda)$ are non-Jacobians.
Maisner~\cite{maisner} extended these ideas to show that for all other values
of $b$ coming from simple isogeny classes, there is a curve with the given
Weil polynomial.

For supersingular surfaces over finite fields of characteristic~$2$,
Question~\ref{Q:main} was answered by Maisner and the second author~\cite{mn2}
by an explicit computation of the zeta functions of all supersingular
curves of genus~$2$, using ideas of van der Geer and van der
Vlugt~\cite{vdGvdV}.
For supersingular surfaces over finite fields of characteristic~$3$,
the question was answered by the first author~\cite{howe:char3}, again
by explicit methods.

McGuire and Voloch~\cite[\S3]{mv} determined which isogeny classes of
split mixed abelian surfaces contain Jacobians, and gave the complete
details of the argument in the case that one factor of the Weil
polynomial of the isogeny class is $x^2\pm2\sqrt{q}x + q$.

In this paper we cover the last steps to get a complete answer
to Question~\ref{Q:main}.
In Part~\ref{Part:split} we deal with the split isogeny classes not
covered by the work of McGuire and Voloch, and in Part~\ref{Part:simple}
we study the simple supersingular isogeny classes that split over~$\ff2$.
These cases are solved with the use of completely different techniques.
For the split case we use a result of Kani~\cite{kani} that characterizes when
two elliptic curves can be tied together along finite subgroups to get a
common covering by a curve of genus two; Kani's result reduces the question
of whether there is a Jacobian isogenous to a product $E\times F$ of
two elliptic curves to the question of whether for some $n>1$ there is an
isomorphism from $E[n]$ to $F[n]$ that is an anti-isometry with respect to
the Weil pairing (and that is `non-reducible', see Section~\ref{S:tying}).
In order to understand the split supersingular case,
we determine the Galois twists of the Dieudonn\'e modules of certain
supersingular elliptic curves.  For the simple supersingular case we use
results of Oort~\cite{oort}, Katsura and Oort~\cite{oort2}, and Ibukiyama,
Katsura, and Oort~\cite{IKO} on supersingular abelian surfaces over the
algebraic closures of finite fields and their polarizations.
Using these results, together with the theory of twists and
work of Hashimoto and Ibukiyama~\cite{hi1} and Ibukiyama~\cite{ibu} on
quaternion hermitian forms, we determine which simple supersingular isogeny
classes contain geometrically non-split principally polarized surfaces.

In our analyses of the supersingular isogeny classes, both split and simple,
it is convenient to assume that the characteristic of the base field is larger
than~$3$.  We may make this assumption because the
characteristic $2$ case is settled in~\cite{mn2}
and the characteristic $3$ case in~\cite{howe:char3}.

\subsection*{Conventions and notation.}
When we speak of a \emph{variety} over a finite field $k$, we mean a variety
defined over the algebraic closure $\kbar$ of $k$ together with
Galois descent data.  By a \emph{morphism} of varieties over a finite
field $k$, we mean a morphism of varieties over $\kbar$ that is
Galois-equivariant.  By a \emph{geometric morphism} of varieties
over $k$, we mean a morphism of varieties over $\kbar$.  Thus, if
$E_1$ and $E_2$ are elliptic curves over a finite field, we will
often speak of geometric isogenies from $E_1$ to $E_2$.
As a consequence of this convention, operators such as
$\Hom$ and $\End$ applied to varieties over $k$
will always refer to $k$-rational homomorphisms and endomorphisms.

If $q$ is a power of a prime $p$, say $q = p^m$, we let $\Q_q$ denote the
unramified degree-$m$ extension of the $p$-adic number $\Q_p$, and we let
$\Z_q$ denote the ring of integers of~$\Q_q$.
We will sometimes denote by $\calA_{(a,b)}$ the isogeny class of
abelian surfaces over $\fq$ with Weil polynomial
$x^4+ax^3+ b x^2 + aq x+q^2$.


\part{Split abelian surfaces as Jacobians}
\label{Part:split}

\section{Introduction}
\label{S:IntroSplit}

In this part of the paper we determine the Weil polynomials of the split
isogeny classes of abelian surfaces that contain Jacobians.

Our first two theorems concern the case of isogeny classes that
contain product surfaces of the form $E_1 \times E_2$, where
$E_1$ and $E_2$ are elliptic curves over $\GF_q$ that are not
isogenous to one another and that are not both supersingular.
Let $s$ and $t$ be the traces
of the Frobenius endomorphisms of $E_1$ and $E_2$,
respectively, so that $s\neq t$ and
so that the Weil polynomial of $E_1\times E_2$
is $$(x^2 - s x + q)(x^2 - t x + q).$$

\begin{theorem}
\label{T:distinct1}
Suppose that neither $s^2$ nor $t^2$ is equal to $4q$
and that $E_1$ and $E_2$ are not both supersingular.
Then there is a Jacobian isogenous to $E_1\times E_2$
if and only if $|s - t| \neq 1$
and $\{q, \{s,t\}\} \neq \{2, \{1,-1\}\}$.
\end{theorem}

\begin{theorem}
\label{T:distinct2}
Suppose that $E_2$ is ordinary and that $s^2 = 4q$,
so that $E_1$ is supersingular.
Then there is a Jacobian isogenous to  $E_1\times E_2$
if and only if $s - t$ is divisible by the square
of an integer greater than~$1$.
\end{theorem}

Theorem~\ref{T:distinct2} was proven by McGuire and Voloch~\cite{mv},
who also mention the special case of Theorem~\ref{T:distinct1} in
which one of the curves is supersingular.  We reprove their
results here for completeness.

Next we  consider isogeny classes that contain squares of
ordinary elliptic curves.

\begin{theorem}
\label{T:squares}
Let $E$ be an ordinary elliptic curve over $\GF_q$ with trace of
Frobenius equal to $t$.  Then there is a Jacobian isogenous to
$E\times E$ if and only if $t^2 - 4 q$ is neither $-3$ nor $-4$
nor $-7$.
\end{theorem}

Finally, we turn to the split supersingular isogeny classes.
We restrict our attention to finite fields of characteristic
greater than $3$, because the characteristic $3$ case is
considered in~\cite{howe:char3} and the
characteristic $2$ case is considered in~\cite{mn2}.
Suppose that $E_1$ and $E_2$ are supersingular elliptic curves over
a finite field $\GF_q$ of characteristic greater than $3$, and
let $s$ and $t$ be the traces of Frobenius of $E_1$ and $E_2$,
respectively.

\begin{theorem}
\label{T:supersingular}
There is a Jacobian isogenous to $E_1\times E_2$ if and only
if $s^2 = t^2$.
\end{theorem}

Our main tool in proving these theorems is a result of Kani~\cite{kani}, which
we review in Section~\ref{S:tying}.
In Section~\ref{S:lemmas} we state and prove a few elementary
lemmas that we will need later in the paper.
In Sections~\ref{S:distinct1},~\ref{S:distinct2}, and~\ref{S:squares}
we prove our theorems for split non-supersingular isogeny classes.
In Section~\ref{S:Dieudonne} we compute certain twists of
the Dieudonn\'e modules of supersingular elliptic curves, and in
Section~\ref{S:supersingular} we use these computations to prove
Theorem~\ref{T:supersingular}.

\section{Tying elliptic curves together along torsion subgroups}
\label{S:tying}

In this section we review a result of Kani that gives necessary and
sufficient conditions for certain split abelian surfaces to be
Jacobians.

Suppose $E_1$ and $E_2$ are elliptic curves over a
field~$k$ and let $n$ be a positive integer.  Suppose
$\psi\colon E_1[n]\to E_2[n]$ is an isomorphism of
group schemes over $k$ that is an anti-isometry with
respect to the Weil pairings on $E_1[n]$ and $E_2[n]$.
Let $A$ be the abelian surface $(E_1\times E_2)/\gra(\psi)$
and let $\varphi\colon E_1\times E_2\to A$ be the natural isogeny.
Then $A$ fits in a diagram
\begin{center}
\leavevmode
\xymatrix{
E_1\times E_2 \ar[d]^{\varphi} \ar[r]^{n}       &  E_1\times E_2  \\
      A                        \ar[r]^{\lambda} & \hat{A}\ar[u]_{\hat{\varphi}}.
}
\end{center}

Here the top arrow is the multiplication-by-$n$ map and
$\hat{A}$ is the dual abelian surface of~$A$.
The existence of the bottom arrow follows from the fact that
the graph of $\psi$ is a maximal isotropic subgroup of the
$n$-torsion of $E_1\times E_2$
(see \cite[Prop.~16.8, p.~135]{milneb}).  In fact, the induced
map $\lambda\colon A\to\hat{A}$ is a polarization, and by looking at the
degrees of the maps in the diagram we see that $\lambda$ is
a principal polarization.  Conversely, if $\lambda$
is a principal polarization of a non-simple abelian surface $A$ over $k$,
then $\lambda$ can be obtained in this way from a pair of
elliptic curves $(E_1,E_2)$ and an anti-isometry $E_1[n]\to E_2[n]$,
for some~$n$.

Kani~\cite{kani} gives a criterion that allows one
to determine when a principally polarized surface $(A,\lambda)$
obtained from an anti-isometry $\psi\colon E_1[n]\to E_2[n]$
is isomorphic to the Jacobian of a curve.  The criterion is
easiest to state when $n$ is a prime.

\begin{theorem}[Kani~{\cite[Thm.~3, p.~95]{kani}}]
\label{T:Kani}
Suppose $n$ is a prime, and let $E_1$, $E_2$, and $\psi$ be
as in the discussion above.
The principally polarized surface
$$(E_1\times E_2)/\gra(\psi)$$
is not a Jacobian if and only if there is an integer $i$
{\rm(}with $0 < i < n${\rm)}
and a geometric isogeny $\varphi\colon E_1\to E_2$
of degree $i(n-i)$ such that
$i \psi = \varphi|_{E_1[n]}.$
\end{theorem}

There is a more complicated criterion when $n$ is composite.
We will only need to use the case $n=4$.

\begin{theorem}[Kani]
\label{T:Kani4}
Suppose $n = 4$ and the characteristic of the base field $k$ is
not equal to~$2$.  Let $E_1$, $E_2$, and $\psi$ be
as in the discussion above.
The principally polarized surface
$$(E_1\times E_2)/\gra(\psi)$$
is not a Jacobian if and only if one of the following
conditions holds{\/\rm:}
\begin{itemize}
\item[\rm(a)] There is a geometric isogeny $\varphi\colon E_1\to E_2$
of degree $3$ such that
$\psi = \varphi|_{E_1[n]}.$
\item[\rm(b)] There are
two order-$2$ subgroups $G_1$ and $G_2$ of $E_1(\GFbar_q)$
and a geometric isomorphism $\varphi\colon E_1\to E_2$
such that the graph of $\psi$ is equal to the set of
points $(x,\varphi(y))$ in $E_1[4](\GFbar_q)\times E_2[4](\GFbar_q)$
such that $x+y\in G_1$ and $x-y\in G_2$.
\end{itemize}
\end{theorem}

\begin{proof}
This follows from Theorem~2.6 of~\cite{kani}.
We make the assumption about the characteristic of $k$ not being~$2$
only so that condition (b) can be stated in terms of groups and not
group-schemes.
\end{proof}

We say that an anti-isometry $\psi\colon E_1[n]\to E_2[n]$ is
\emph{reducible}  if $(A,\lambda)$ is
not a Jacobian.  If $\psi$ is not reducible, then
we refer to the process of constructing a Jacobian from
$E_1$, $E_2$, and $\psi$
as \emph{tying $E_1$ and $E_2$ together along their
$n$-torsion subgroups via $\psi$.}

\section{Useful lemmas}
\label{S:lemmas}

In this section we present some lemmas that will be helpful in later
sections.

Suppose that $E$ is an elliptic curve over $\GF_q$ with trace of
Frobenius $t$, and suppose that $t^2\neq 4q$.  Let $\pi$ be the
Frobenius endomorphism of $E$ and let $R$ be the subring $\Z[\pi]$
of $\End(E)$.  Then $R$ is an imaginary quadratic order of
discriminant $t^2 - q$.  Let $\calO$ be
the integral closure of~$R$ in $R\otimes\Q$.  The endomorphism ring of $E$
is an order that is contained in $\calO$ and that contains~$R$.  Let
$\ell$ be a prime integer.  We say that $E$ is \emph{maximal}
at $\ell$ if $\ell$ does not divide the index of $\End(E)$ in $\calO$.
We say that $E$ is \emph{minimal} at $\ell$ if $\ell$ does not
divide the index of $R$ in $\End(E)$.

Given an $E$ as above, let $\ell$ be a prime that does not divide
both $t$ and $q$.  Then it follows from~\cite[Thm.~4.2]{wa}
that there is a curve isogenous to $E$ that is minimal at~$\ell$.

\begin{lemma}
\label{L:MinimalPolynomials}
Let $E$ be an elliptic curve over $\GF_q$ whose Weil polynomial
is $x^2 - tx + q$, and suppose that $t^2 \neq 4q$.
Let $\ell$ be a prime that does not divide both $t$ and $q$,
and suppose that $E$ is minimal at~$\ell$.
Then the minimal polynomial of Frobenius acting on
$E[\ell]$ is $x^2 - tx + q\in\GF_\ell[x]$.
\end{lemma}

\begin{proof}
The characteristic polynomial of Frobenius on $E[\ell]$ is
$x^2 - tx + q$, and the only way that this might not be the
minimal polynomial of Frobenius is if $t^2 - 4q \equiv 0 \bmod \ell$.

Suppose $t^2 - 4q \equiv 0 \bmod \ell$.
Let $b$ be an integer such that $t \equiv 2b \bmod \ell$,
and let $\pi$ be the Frobenius endomorphism of~$E$.
Then the characteristic polynomial of
$\pi$ on $E[\ell]$ is $(x - b)^2$, and the minimal polynomial
will be $x - b$ if and only if $(\pi - b)/\ell$ lies in $\End(E)$.
But the index of $\End(E)$ in $\Z[\pi]$ is coprime to $\ell$
by assumption.
\end{proof}

\begin{lemma}
\label{L:QuadraticPrimes}
Let $\ell$ be either $4$ or an odd prime and let $K$
be an imaginary quadratic field whose discriminant is not equal
to $-\ell$.  Then there are infinitely many rational
primes $m$ that split in $K$ and that are nonsquares
modulo~$\ell$.
\end{lemma}

\begin{proof}
Let $\Delta$ be the discriminant of $K$. The splitting of a
rational prime in $K$ depends only on its congruence class
modulo $\Delta$.  If $\Delta$ is coprime to $\ell$, then
we can choose a congruence class modulo $\ell\Delta$ such
that every prime $m$ in this congruence class splits in
$K$ and is a nonsquare modulo $\ell$.
So from this point on we consider the case where $\ell$ and $\Delta$
have a common factor.

We first consider the case where $\ell$ is prime.

Write $\Delta = -2^e \ell D$ for some odd positive $D$, and suppose
that $D>1$.  If $m$ is a prime that is congruent to $1$ modulo $8$,
that is not a square modulo $\ell$, and such that the Jacobi symbol
$(m/D)$ is $-1$, then we have
$$\left(\frac{\Delta}{m}\right) =
\left(\frac{-1}{m}\right)
\left(\frac{2}{m}\right)^e
\left(\frac{\ell}{m}\right)
\left(\frac{D}{m}\right)
=
1\cdot 1\cdot
\left(\frac{m}{\ell}\right)
\left(\frac{m}{D}\right)
=
1.$$
So when $D>1$, there are infinitely many primes that split in
$K$ and that are not squares modulo~$\ell$.

Suppose $D=1$, so that either $\Delta = -8\ell$ or
$\Delta = -4\ell$.  Suppose $\Delta = -8\ell$.
Then if $m$ is a prime that is congruent to $5$ modulo $8$ and
that is a nonsquare modulo $\ell$, then
$$\left(\frac{\Delta}{m}\right) =
\left(\frac{-1}{m}\right)
\left(\frac{2}{m}\right)
\left(\frac{\ell}{m}\right)
=
1\cdot (-1)\cdot
\left(\frac{m}{\ell}\right)
=
1.$$
On the other hand, suppose that $\Delta = -4\ell$.
This can only happen if $\ell\equiv 1\bmod 4$.
If $m$ is a prime that is congruent to $3$ modulo $4$ and
that is a nonsquare modulo $\ell$, then
$$\left(\frac{\Delta}{m}\right) =
\left(\frac{-1}{m}\right)
\left(\frac{\ell}{m}\right)
=
(-1)\cdot
\left(\frac{m}{\ell}\right)
=
1.$$
Thus, in every case there are infinitely many primes that split in
$K$ and that are not squares modulo~$\ell$.

Next we consider the case where $\ell = 4$.

Suppose $\Delta = -4D$ for some odd $D > 1$.
This can only be the case if $D\equiv 1 \bmod 4.$
If $m$ is a prime that is congruent to $3$ modulo $4$
and such that the Jacobi symbol
$(m/D)$ is $-1$, then we have
$$\left(\frac{\Delta}{m}\right) =
\left(\frac{-1}{m}\right)
\left(\frac{D}{m}\right)
=
(-1)\cdot
\left(\frac{m}{D}\right)
=
1,$$
and we see that there are
infinitely many primes that split in
$K$ and that are not squares modulo~$\ell$.

Suppose $\Delta = -8D$ for some odd $D>0$ (with $D=1$ being allowed).
Suppose $D\equiv 1 \bmod 4$.  Then if $m$ is a prime that is $3$
modulo $8$ and such that
the Jacobi symbol $(m/D)$ is $1$,  we have
$$\left(\frac{\Delta}{m}\right) =
\left(\frac{-1}{m}\right)
\left(\frac{2}{m}\right)
\left(\frac{D}{m}\right)
=
(-1)\cdot (-1) \cdot
\left(\frac{m}{D}\right)
=
1.$$
On the other hand, if $D\equiv 3 \bmod 4$,
we can consider primes $m$ that are $7$ modulo $8$ and such that
the Jacobi symbol $(m/D)$ is $-1$, so that
$$\left(\frac{\Delta}{m}\right) =
\left(\frac{-1}{m}\right)
\left(\frac{2}{m}\right)
\left(\frac{D}{m}\right)
=
(-1)\cdot 1 \cdot
\left(\frac{m}{D}\right)
=
1.$$
Again we see that in both cases there are
infinitely many primes that split in
$K$ and that are not squares modulo~$\ell$.
\end{proof}

Next we present a lemma that provides us with a large
supply of anti-isometries.
This lemma and the one that follows it
will refer to integers $\ell$ that are
assumed to be either prime or equal to~$4$.
It will be convenient
to define $\ell^*$ to be the unique prime divisor of such an
integer $\ell$.

\begin{lemma}
\label{L:PreAntiIsometry}
Let $E_1$ and $E_2$ be elliptic curves over a finite field
$\GF_q$ and let $s$ and $t$ be their traces of Frobenius.
Suppose that $|s-t|$ is neither $0$ nor $1$ and that neither $s^2$ nor $t^2$
is equal to~$4q$.
Let $\ell$ be a divisor of $s - t$ that is either $4$ or a prime,
and assume that $\ell$ is coprime to $q$ if either $E_1$ or $E_2$
is supersingular.
If $\ell=2$ then let $m=1${\rm;} otherwise,
let $m$ be a positive integer coprime to $\ell$
whose image in $(\Z/\ell\Z)$ is a nonsquare.
Suppose that $E_1$ and $E_2$ are minimal at $\ell^*$, and that
$E_1'$ is an elliptic curve that is $m$-isogenous to $E_1$.
Then either there is an $\GF_q$-defined anti-isometry
from $E_1[\ell]$ to $E_2[\ell]$, or there is one from
$E_1'[\ell]$ to $E_2[\ell]$.
\end{lemma}

\begin{proof}
By Lemma~\ref{L:MinimalPolynomials}
the minimal polynomials of Frobenius on $E_1[\ell^*]$
and $E_2[\ell^*]$ are both equal to $x^2 - t x + q\in\GF_\ell[x]$.
It follows easily that there are
points $P_1\in E_1[\ell](\GFbar_q)$ and $P_2\in E_2[\ell](\GFbar_q)$
that generate the Galois modules
$E_1[\ell](\GFbar_q)$ and $E_2[\ell](\GFbar_q)$, respectively.

We claim that the $k$-group schemes $E_1[\ell]$ and $E_2[\ell]$
are isomorphic.  If $p\neq\ell^*$
then we can see this by defining an isomorphism
$E_1[\ell] \to E_2[\ell]$ by sending
$P_1$ to $P_2$ and extending by Galois equivariance.
If $p=\ell^*$, then $E_1[\ell]$ and $E_2[\ell]$ are both products of a
reduced group scheme of rank $\ell$ and a local group scheme of
rank $\ell$.  On each of the reduced group schemes, the Frobenius acts
as multiplication-by-$t$, so the reduced subschemes of $E_1[\ell]$
and $E_2[\ell]$ are isomorphic.  But the local subschemes are the
duals of the reduced subschemes, so the local subschemes are
isomorphic as well.  Thus $E_1[\ell]$ and $E_2[\ell]$
are isomorphic.

Let $\psi\colon E_1[\ell]\to E_2[\ell]$ be an isomorphism of
group schemes over $\GF_q$.
For each $i=1,2$ let $e_i$ be the Weil pairing
$E_i[\ell]\times E_i[\ell]\to \Bmu_\ell$.
Then there is an element $r$ of
$\Aut \Bmu_\ell \cong (\Z/\ell\Z)^*$ such
that the diagram commutes:
\begin{center}
\leavevmode
\xymatrix@C=2cm{
E_1[\ell]\times E_1[\ell] \ar[r]^{\quad (\psi,\psi)\quad } \ar[d]^{e_1} &      E_2[\ell]\times E_2[\ell]\ar[d]^{e_2} \\
\Bmu_\ell\ar[r]^{r}                                                     & \Bmu_\ell
}
\end{center}

Let $F$ be an elliptic curve isogenous to $E_1$, and
suppose $\varphi\colon F\to E_1$ is an isogeny of degree coprime
to $\ell$.  If we replace $E_1$ by $F$ and $\psi$ by
$\psi\circ\varphi$,
then $r$ is replaced with $r$ times the degree of $\varphi$.
Using multiplication by integers in $\End(E_1)$, we can modify
$r$ in this way by arbitrary squares in $(\Z/\ell\Z)^*$.
Using the isogeny $E_1'\to E_1$, we can modify
$r$ by the integer~$m$, which is a nonsquare
modulo $\ell$ when $\ell\neq 2$.  Therefore we can modify $r$ so
that it is equal to $-1$; in other words, we can find an
anti-isometry that maps either $E_1[\ell]$ or $E_1'[\ell]$
to $E_2[\ell]$.
\end{proof}

The next lemma is a useful special case of
Lemma~\ref{L:PreAntiIsometry}.

\begin{lemma}
\label{L:AntiIsometry}
Let $E_1$ and $E_2$ be elliptic curves over a finite field
$\GF_q$ and let $s$ and $t$ be their traces of Frobenius.
Suppose that $|s-t|$ is neither $0$ not $1$ and that neither
$s^2$ nor $t^2$ is equal to~$4q$.
Write $s^2 - 4q = f_1^2\Delta_1$ and $t^2 - 4q = f_2^2 \Delta_2$
for integers $f_i$ and fundamental discriminants~$\Delta_i$.
Let $\ell$ be a divisor of $s - t$ that is either $4$ or a prime,
and assume that $\ell$ is coprime to $q$ if either $E_1$ or $E_2$
is supersingular.
Suppose that $\Delta_1$ and $\Delta_2$ are not both equal
to~$-\ell$.
Then there are elliptic curves $F_1$ and $F_2$, isogenous
to $E_1$ and $E_2$, respectively, for which there is
an $\GF_q$-defined anti-isometry $F_1[\ell]\to F_2[\ell]$.
\end{lemma}

\begin{proof}
By symmetry, we may assume that $\Delta_1\neq-\ell$.

Replace $E_1$ and $E_2$ with isogenous curves that
are minimal at~$\ell^*$.  If $\ell=2$ then Lemma~\ref{L:PreAntiIsometry}
shows that there is an anti-isometry $E_1[\ell]\to E_2[\ell]$,
so we may assume that $\ell>2$.

Let $m$ be a prime number.  If $m$ splits in
$\End(E_1)\otimes\Q\cong\Q(\sqrt{\Delta_1})$ then there
is a degree-$m$ isogeny from $E_1$ to some other elliptic
curve $F$.  (Indeed, Ito~\cite{ito}
proves that there is a degree-$m$ isogeny from $E_1$
to some $F$ if and only if either $m$ splits or ramifies
in $\End(E_1)\otimes\Q\cong\Q(\sqrt{\Delta_1})$ or $m$ divides~$f_1$.)
Lemma~\ref{L:QuadraticPrimes} says that there are
infinitely many $m$ that are nonsquares modulo $\ell$ and
that split in $\End(E_1)\otimes\Q$, so we
know there is an elliptic curve $E_1'$ that is $m$-isogenous
to $E_1$ for some $m$ that is not a square modulo~$\ell$.
Using Lemma~\ref{L:PreAntiIsometry}, we see that there is either
an anti-isometry $E_1[\ell]\to E_2[\ell]$ or an anti-isometry
$E_1'[\ell]\to E_2[\ell]$, and we are done.
\end{proof}

We end with a lemma that will help us show that certain
anti-isometries are not reducible.

\begin{lemma}
\label{L:twists}
Suppose $E$ and $F$ are ordinary elliptic curves over a finite field $k$,
let $F'/k$ be a twist of $F$, and let $\chi\colon F'\to F$ be a
geometric isomorphism.  Suppose that $E$ and $F'$ are isogenous over~$k$.
Let $\ell$ be a prime, and suppose $\psi\colon E[\ell]\to F[\ell]$ is a
Galois-equivariant anti-isometry.  If $\psi$ is reducible, then
$\chi|_{F'[\ell]}$ is Galois equivariant.
\end{lemma}

\begin{proof}
Suppose $\psi$ is reducible.  Then Theorem~\ref{T:Kani} shows
that there is an integer $i$ and a
geometric isogeny $\varphi\colon E\to F$ of degree $i(n-i)$
such that
$i \psi = \varphi|_{E[\ell]}.$
The left-hand side of this equality is Galois equivariant,
so $\varphi|_{E[\ell]}$ is Galois equivariant.

Every geometric isogeny $E\to F$ can be written as the composition
of a geometric isogeny $E\to F'$ with the geometric isomorphism
$\chi\colon F'\to F$, so we may write
$\varphi = \chi\circ\varphi'$ for a geometric isogeny
$\varphi'\colon E\to F'$.  Since $E$ and $F'$ are ordinary,
all of their endomorphisms are defined over $k$,
and all isogenies from $E$ to $F'$ are defined over $k$.
It follows that
$\varphi'$ is Galois equivariant.  Also, $\varphi'$
gives an isomorphism $E[\ell]\to F'[\ell]$, so
$\chi$ induces a Galois equivariant
isomorphism from $F'[\ell]$ to $F[\ell]$.
\end{proof}

\section{Proof of Theorem \ref{T:distinct1}}
\label{S:distinct1}

Suppose that $|s - t| = 1$.  Then a result of Serre
(see~\cite[Lem.~1]{lauter} or~\cite[Thm.~1(a)]{hola})
shows that there is no Jacobian isogenous to $E_1\times E_2$.

Suppose that $|s - t| > 1$.
The remainder of the
proof of Theorem~\ref{T:distinct1} breaks into five cases:
\begin{enumerate}
\item
$E_1$ and $E_2$ are geometrically non-isogenous.
\item
$E_1$ and $E_2$ become isogenous to one another over a degree-$2$
extension.
\item
$E_1$ and $E_2$ become isogenous to one another over a degree-$3$
extension.
\item
$E_1$ and $E_2$ become isogenous to one another over a degree-$4$ extension,
   but not over a degree-$2$ extension.
\item
$E_1$ and $E_2$ become isogenous to one another over a degree-$6$ extension,
   but not over a degree-$2$ or degree-$3$ extension.
\end{enumerate}

To see that these cases include all possibilities, we note that
if $E_1$ and $E_2$ are geometrically isogenous to one another
then they must both be ordinary (because the hypotheses of the
theorem preclude them from both being supersingular).  Let $\pi_1$
and $\pi_2$ be the Weil numbers of $E_1$ and $E_2$,
considered as elements of $\Qbar$.  If $E_1$ and $E_2$
become isogenous over a degree-$n$ extension (and no smaller
extension), then
$\pi_1^n$ and $\pi_2^n$ are conjugate quadratic integers;
replacing $\pi_1$ with its conjugate, if necessary, we
may assume that $\pi_1^n = \pi_2^n$.  Then the two quadratic
fields $\Q(\pi_1)$ and $\Q(\pi_2)$ are equal, and they
contain the primitive $n$-th root of unity $\pi_1/\pi_2$.
This restricts the possible values for $n$ to be $2$, $3$, $4$, and~$6$.

We will consider these cases separately.  In each case, we
will denote the characteristic of $\GF_q$ by $p$.

\subsection{Case 1: $E_1$ and $E_2$ are geometrically non-isogenous.}

Pick a prime $\ell$ dividing $s-t$.  Since $E_1$ and $E_2$
are not both supersingular, $\ell$ is not equal to $p$
if either curve is supersingular.

Let $\pi_1$ and $\pi_2$ be the Weil numbers of $E_1$ and
$E_2$, respectively, and let $K_1$ and $K_2$ be the
imaginary quadratic fields generated by $\pi_1$ and $\pi_2$.
If one of these fields has discriminant unequal to $-\ell$,
then Lemma~\ref{L:AntiIsometry} shows that we can replace
$E_1$ and $E_2$ with isogenous curves for which there
is an $\GF_q$-defined anti-isometry $E_1[\ell]\to E_2[\ell]$.
Since $E_1$ and $E_2$ are geometrically non-isogenous
by assumption, Theorem~\ref{T:Kani} shows that we can
tie $E_1$ and $E_2$ together along their $\ell$-torsion
to get a genus-$2$ curve with Jacobian isogenous to $E_1\times E_2$.

We are left to consider the case where $K_1$ and $K_2$ are
both isomorphic to the imaginary quadratic field $K$ of
discriminant~$-\ell$.
In this case, we may view $\pi_1$ and $\pi_2$ as elements of~$K$.

If $E_1$ and $E_2$ are both ordinary, then their Weil numbers
must differ from one another by a root of unity (and perhaps
complex conjugation).  But then $E_1$ and $E_2$ become isogenous
to one another after a base field extension, contradicting
our hypotheses.

So suppose one of our elliptic curves, say $E_1$, is supersingular and the
other is ordinary.  We have already noted that in this case $\ell\neq p$.
Now, we know the possible Weil numbers for supersingular elliptic
curves --- see~\cite{schoof} or~\cite[Thm.~4.2]{wa} for example --- and the
only way a supersingular Weil number can generate an imaginary quadratic field
of prime discriminant unequal to $-p$ is if that field is $\Q(\sqrt{-3})$ and
if $p\neq 1\bmod 3$.  But then $p$ does not split in $\Q(\sqrt{-3})$, so
there are no ordinary elliptic curves over $\GF_q$ with CM by $\Q(\sqrt{-3})$,
contradicting the existence of $E_2$.

\subsection{Case 2: $E_1$ and $E_2$ become isogenous to one another
over a degree-$2$ extension.}

In this case, there is an integer $t$ such that the Weil polynomials
of $E_1$ and $E_2$ are $x^2 - t x + q$ and $x^2 + t x + q$.
Also, $E_1$ and $E_2$ are both ordinary.  Let $\Delta = t^2 - 4q$
and let $R$ be the imaginary quadratic order of
discriminant~$\Delta$.

We have two arguments that each cover many cases, and that
together cover all but one case.

\begin{lemma}
\label{L:ClassNumber}
If the class number of $R$ is greater than $1$,
there is a Jacobian isogenous to $E_1\times E_2$.
\end{lemma}

\begin{proof}
If the class number of $R$ is greater than~$1$, we can choose
$E_1$ and $E_2$ to have endomorphism ring $R$, and to
be geometrically non-isomorphic.
Then $E_1$ and $E_2$ are minimal at~$2$, so
Lemma~\ref{L:PreAntiIsometry} shows that there is
an anti-isometry $\psi\colon E_1[2]\to E_2[2]$.
Theorem~\ref{T:Kani} shows that $\psi$ is not
reducible because there are no geometric isomorphisms $E_1\to E_2$.
Thus we can
tie $E_1$ and $E_2$ together along their $2$-torsion to get a
Jacobian isogenous to $E_1\times E_2$.
\end{proof}

\begin{lemma}
\label{L:Trace}
If $|t|>1$ then there is a Jacobian isogenous to $E_1\times E_2$.
\end{lemma}

\begin{proof}
Let $\ell$ be a prime divisor of $t$.

First suppose that $\ell$ is odd.
Note that $\ell$ is not equal to~$p$,
because $\ell$ divides $t$ and $E_1$ and $E_2$ are not both
supersingular.  It follows that $\ell$ does not divide
$\Delta = t^2 - 4q$.  Then
Lemma~\ref{L:AntiIsometry} shows that we can replace
$E_1$ and $E_2$ with isogenous curves for which there
is an $\GF_q$-defined anti-isometry $\psi\colon E_1[\ell]\to E_2[\ell]$.
Let $E_2'$ be the quadratic twist of $E_2$ and let
$\chi\colon E_2'\to E_2$ be the standard geometric isomorphism.
Then Lemma~\ref{L:twists} shows that if $\psi$ is reducible,
then $\chi$ induces a Galois-equivariant isomorphism
from $E_2'[\ell]$ to $E_2[\ell]$.  But from the
definition of the quadratic twist, we know that
$\chi(P^\sigma) = - (\chi(P))^\sigma$ for all geometric
points $P$, where $\sigma$ denotes the $q$-th power
Frobenius automorphism of $\GFbar_q$.  From this it
is clear that $\chi$ does \emph{not} give a
Galois-equivariant isomorphism $E_2'[\ell]\to E_2[\ell]$,
because $\ell > 2$.
Thus $\psi$ is not reducible, and we can
tie $E_1$ and $E_2$ together along their $\ell$-torsion
via $\psi$ to get a
Jacobian isogenous to $E_1\times E_2$.

Next suppose $\ell = 2.$  Lemma~\ref{L:AntiIsometry} shows that
we can replace $E_1$ and $E_2$ with isogenous curves for which
there is an anti-isometry $\psi\colon E_1[4]\to E_2[4]$.
According to Lemma~\ref{L:twists}, there are two ways in which
$\psi$ might be reducible.
In the first way, there is a geometric isogeny
$\varphi\colon E_1\to E_2$ of degree $3$ such that
$\varphi|_{E_1[\ell]} = \pm\psi$.
But as in the argument for odd~$\ell$,
we obtain a contradiction from the facts that
$\pm\psi$ is Galois equivariant
while $\varphi|_{E_1[\ell]}$ is not.

The other way that $\psi$
can be reducible is if there is a geometric isomorphism
$\varphi\colon E_1\to E_2$ and two order-$2$ subgroups
$G_1$ and $G_2$ of $E_1(\GFbar_q)$ such that the
graph of $\psi$ is equal to the
set of $(x,\varphi(y))$ in $E_1[4](\GFbar_q)\times E_2[4](\GFbar_q)$
such that $x+y\in G_1$ and $x-y\in G_2$.

In particular, if $E_2$ is not the quadratic twist of $E_1$ then $\psi$
is not reducible.  Suppose $E_2$ is the quadratic twist of $E_1$ and
that $\psi$ is reducible.  Identify $E_2[4]$ with $E_1[4]$ provided with
the negative Galois action.  Pick $4$-torsion points $X$ and $Y$ of
$E_1(\GFbar_q)$ such that $G_1 = \langle 2X\rangle$ and
$G_2 = \langle 2Y\rangle$.  Then $\psi(aX+bY) = \varphi(aX - bY)$
for some automorphism $\varphi$ of $E_2$.  Since $E_2$
is ordinary, all of its automorphisms are
defined over $\GF_q$, so $\psi$ is Galois equivariant if and only if
$\varphi^{-1}\circ \psi$ is Galois equivariant, and one can
show that this is the case
if and only if $2X^\sigma = 2Y$ and $2Y^\sigma = 2X$;
here $\sigma$ denotes the Frobenius
element of $\Gal(\GFbar_q/\GF_q)$.

So assume that $2X^\sigma = 2Y$ and $2Y^\sigma = 2X$.
Under this assumption, we can compute the number of
reducible Galois equivariant anti-isometries.
The only choices
for $G_1$ and $G_2$ are $G_1=\langle 2X\rangle$,
$G_2=\langle 2Y\rangle$ and
$G_1=\langle 2Y\rangle$, $G_2=\langle 2X\rangle$,
and swapping $G_1$ and $G_2$ is equivalent to replacing $\varphi$
with $-\varphi$.
Thus we see that
the automorphism group of $E_2$ acts transitively on the set
of reducible Galois-equivariant anti-isometries.

Let us compute the number of Galois equivariant anti-isometries
there are from $E_1[4]$ to $E_2[4]$. To make the computation
simpler, we replace $Y$ with $X^\sigma$ (this does not change
the group $G_2$), and write
$Y^\sigma = \eps X + f Y$, where $\eps = \pm 1$ and
$f\in\{0,2\}$.  Note that it follows that
the characteristic polynomial
of Frobenius on the $4$-torsion of $E_1$ is
$x^2 + fx - \eps$.  Since the Weil polynomial of $E_1$
is $x^2 - tx + q$, we see that $\eps\equiv -q\bmod 4$ and
$f\equiv t \bmod 4$.

A Galois-equivariant map $\chi\colon E_1[4]\to E_2[4]$ that sends $X$
to $aX+bY$ must send $Y$ to $bq X - (a+tb)Y$,
and $\chi$ will be an anti-isometry
if and only if $a^2 + tab + q b^2 \equiv 1\bmod 4$.
There are four pairs $(a,b)$ of elements of $(\Z/4\Z)$ that satisfy
this condition if $q \equiv -1\bmod 4$ and eight if
$q\equiv 1\bmod 4$, so there are
either exactly $4$ or exactly $8$ Galois-equivariant anti-isometries.

If $\#\Aut E_2 = 2$ then we have only $2$ reducible
Galois-equivariant anti-isometries, so there
are at least two nonreducible ones.  If $\#\Aut E_2 = 6$ then either there
are $8$ Galois-equivariant anti-isometries in total and we can choose
a nonreducible one, or else there are only $4$ Galois-equivariant
anti-isometries, in which case $\Aut E_2$ cannot act faithfully on the
reducible Galois-equivariant anti-isometries; but since the automorphism
$-1$ doesn't act trivially, the kernel of the action must be of
order $3$, and once again we see that
there are only two reducible Galois-equivariant
anti-isometries.

If $\#\Aut E_2 = 4$ then $\Delta = t^2 - 4 q = -4$, and
the Frobenius can be written $\pi = (t/2) + i$.
Note that therefore $q\equiv 1 \bmod 4$. We noted
above that in this case there are eight Galois-equivariant
anti-isometries,
so we have four nonreducible ones to choose from.
\end{proof}

The only situations not covered by Lemmas~\ref{L:ClassNumber}
and~\ref{L:Trace} are those in which $t=1$ and in which
$\Delta$ lies in the set
$$\{-3,-12,-27,-4,-16,-7,-28,-8,-11,-19,-43,-67,-163\}.$$
That means the only cases left to consider are
those in which $(q,\Delta)$ is one of
$$\{(2,-7),(3,-11),(5,-19),(7,-27),(11,-43),(17,-67),(41,-163)\}.$$

Suppose $q$ is odd.  Then characteristic polynomial
of Frobenius is congruent to $x^2 + x + 1$ modulo $2$.
In this case there are three Galois-equivariant anti-isometries
from $E_1[2]$ to $E_2[2]$.  If they are all reducible then we must have
$\#\Aut E_1 = 6$, so that $\Delta = -3$.
But this is not one of the $\Delta$'s on our list of $(q,\Delta)$
pairs.

Finally we are left with $t = 1$ and $q = 2$.
We find, by explicitly enumerating the genus-$2$ curves
over $\GF_2$, that none of them has Weil polynomial
$(x^2 + x + 2) (x^2 - x + 2).$

\subsection{Case 3: $E_1$ and $E_2$ become isogenous to one another
over a degree-$3$ extension.}

In this case the Weil numbers of $E_1$ and $E_2$ must both
live in $\Q(\omega)$, where $\omega^2 + \omega + 1 = 0$,
and they can be chosen so that they differ multiplicatively by $\omega$.
So let us write
\begin{align*}
\pi_1 &= a + b\omegabar\\
\pi_2 &= \omega\pi_1 = b + a\omega.
\end{align*}
Note that
\begin{align*}
q         &= \pi_1\mybar{\pi_1} = a^2 - ab + b^2\\
s         &= 2a - b\\
t         &= 2b - a\\
\Delta_1  &= -3b^2\\
\Delta_2  &= -3a^2\\
s - t     &= 3(a-b).
\end{align*}
We observe several facts.
First, we see that $\Delta_1$ and $\Delta_2$ cannot be equal;
if they were equal, then we would have $b=-a$,
and $3$ would  divide both $s$ and $q$, contradicting
the ordinariness of $E_1$.  Second, we note that the
same reasoning shows that $a$ and $b$ are coprime to each other.
Third, we see that $q$ is congruent to $1$ modulo $3$, and
in particular, the characteristic of $\GF_q$ is not~$3$.
And fourth, we see that $q$ is odd, so the
characteristic of $\GF_q$ is not~$2$.

Suppose $a$ and $b$ are both odd, so that $2$ divides $s-t$.
Replace $E_1$ and $E_2$ with isogenous curves whose endomorphism
rings have discriminants $\Delta_1$ and $\Delta_2$, respectively,
so that in particular $E_1$ and $E_2$ are geometrically
non-isomorphic.  Then there is an anti-isometry
$\psi\colon E_1[2]\to E_2[2]$, and since $E_1$ and $E_2$
are geometrically non-isomorphic, Theorem~\ref{T:Kani}
shows that $\psi$ is not reducible.  Thus we may tie $E_1$
and $E_2$ together along their $2$-torsion to get a genus-$2$
curve.

We are left to consider the case in which one of $a$ and $b$ is even.
By symmetry, we may assume that $b$ is even.

Suppose that $a$ is not a multiple of~$3$.
Replace $E_2$ with an isogenous curve that has complex
multiplication by $\Z[\omega]$.  Since $a$ is not a multiple
of~$3$, we see that $E_2$ is minimal at~$3$.
Applying Lemma~\ref{L:AntiIsometry}, we find that we can
replace $E_1$ by an isogenous curve so that there is
an anti-isometry $\psi\colon E_1[3]\to E_2[3]$.

Let $F$ be the cubic twist of $E_2$ that is isogenous to $E_1$,
and let $\chi\colon F\to E_2$ be a geometric isomorphism.
Lemma~\ref{L:twists} shows that if $\psi$ is reducible, then
$\chi$ induces an isomorphism from $F[3]$ to $E_2[3]$ as
group schemes over $\GF_q$.

We know that $E_2$ can be written in the form
$y^2 = x^3 + e$ for some $e\in\GF_q$,
and the twist $F$ of $E_2$ can be written
$y^2 = cx^3 + e$ for some $c\in\GF_q$ that is not a cube.
Then the geometric isomorphism $\chi$ can be
taken to be $(x,y) \mapsto (dx,y)$, where $d\in\GFbar_q$
satisfies $d^3 = c$.  But then it is clear that $\chi$ will
not induce a Galois-equivariant isomorphism $F[3]\to E_2[3]$
if $F[3]$ contains an element with nonzero $x$-co\"ordinate.
Since $F[3]$ clearly contains such an element,
$\psi$ must not be reducible, so we can tie
$E_1$ and $E_2$ together along their $3$-torsion.

Finally, suppose that $a$ is divisible by~$3$.
Replace $E_2$ with an isogenous elliptic curve whose
endomorphism ring has discriminant~$\Delta_2$, and
replace $E_1$ with an isogenous elliptic curve
with complex multiplication by $\Z[\omega]$.
Since $b$ is even, there is a $2$-isogeny
from $E_1$ to an elliptic curve whose endomorphism
ring is $\Z[\sqrt{-3}]$.  Lemma~\ref{L:PreAntiIsometry}
shows that there is an anti-isometry $\psi$ from either $E_1[3]$ or
$E_1'[3]$ to $E_2[3]$.  Lemma~\ref{L:twists} shows that if
this isometry is reducible, there must be a geometric
$2$-isogeny from $E_1$ or $E_1'$ to $E_2$.  But by looking
at the discriminants of the endomorphism rings of these curves,
we see that every isogeny from $E_1$ or $E_1'$ to $E_2$ must
have degree divisible by~$3$.  Thus $\psi$ is not reducible,
so we may use it to produce a genus-$2$ curve whose Jacobian
is isogenous to~$E_1\times E_2$.

\subsection{Case 4: $E_1$ and $E_2$ become isogenous to one another
over a degree-$4$ extension, but not over a degree-$2$ extension.}

In this case the Weil numbers of $E_1$ and $E_2$
must live in $\Q(i)$, where $i^2 = -1$, and they may be chosen
so that they differ multiplicatively by~$i$.
So let us write
\begin{align*}
\pi_1 &=  a + bi\\
\pi_2 &= -b + ai\\
\intertext{so that we have}
q        &=  a^2 + b^2\\
s        &=  2a\\
t        &= -2b\\
\Delta_1 &= -4b^2\\
\Delta_2 &= -4a^2.
\end{align*}

If $b$ were equal to $\pm a$ then $q$ and $s$ would both be even,
contradicting our assumption that $E_1$ is ordinary.
Therefore $\Delta_1$ and $\Delta_2$ are not equal to one another,
so if we pick $E_1$ and $E_2$ with minimal endomorphism rings,
they will be geometrically non-isomorphic.  Since $s-t$
is even, we can tie $E_1$ and $E_2$ together
along their $2$-torsion.

\subsection{Case 5: $E_1$ and $E_2$ become isogenous to one another
over a degree-$6$ extension,  but not over a degree-$2$ or degree-$3$
extension.}

In this case the two Weil numbers must live in $\Q(\omega)$,
where $\omega^2 + \omega + 1 = 0$,  and they can be chosen
to differ multiplicatively by a primitive sixth root of unity,
such as $-\omega$.  So let us write
\begin{align*}
\pi_1 &= a + b\omegabar\\
\pi_2 &= -\omega\pi_1 = -b - a\omega\\
\intertext{so that we have}
q        &= \pi_1\mybar{\pi_1} = a^2 - ab + b^2\\
s        &= 2a - b\\
t        &= a - 2b\\
\Delta_1 &= -3b^2\\
\Delta_2 &= -3a^2\\
s-t      &= a+b.
\end{align*}
As in Case 3, we see that $\Delta_1\neq \Delta_2$,
that $(a,b) = 1$, and that the characteristic of
$\GF_q$ is not~$2$.  Since $(a,b)=1$, at least one of
$a$ and $b$ is odd; by symmetry, we may assume that $a$ is odd.
Let $\ell$ be the smallest
prime divisor of $s-t = a+b$.  Note that neither $a$ nor $b$
can be divisible by $\ell$, so both $E_1$ and $E_2$ are
automatically minimal at~$\ell$.

Let us replace $E_2$ with an isogenous curve that has complex
multiplication by $\Z[\omega]$.
We will show that we can replace $E_1$ with an isogenous
curve for which there is
an anti-isometry $E_1[\ell]\to E_2[\ell]$.

First suppose that $\ell\neq 3$.
Since $E_1$ has complex multiplication by an order in
$\Q(\omega)$ and since $\ell\neq3$, we know from
Lemma~\ref{L:QuadraticPrimes} that there is a prime $m\equiv 2\bmod 3$
for which there is an elliptic curve $E_1'$ that is $m$-isogenous
to $E_1$.  Applying Lemma~\ref{L:PreAntiIsometry}, we find that
there is an anti-isometry from either $E_1[\ell]$ or
$E_1'[\ell]$ to $E_2[\ell]$.

On the other hand, suppose that $\ell=3$.
Since $\ell$ was chosen to be the smallest prime divisor of
$s-t = a+b$, we know that $a+b$ is odd; since $a$ is
odd, we know that $b$ is even.  It follows that $E_1$
is $2$-isogenous to some other curve $E_1'$;
applying Lemma~\ref{L:PreAntiIsometry}, we find that
there is an anti-isometry from either $E_1[\ell]$ or
$E_1'[\ell]$ to $E_2[\ell]$.

Let $F$ be the sextic twist of $E_2$ that is isogenous to $E_1$,
and let $\chi\colon F\to E_2$ be a geometric isomorphism.
Lemma~\ref{L:twists} shows that if $\psi$ is reducible, then
$\chi$ induces an isomorphism from $F[\ell]$ to $E_2[\ell]$ as
group schemes over $\GF_q$.

We know that $E_2$ can be written in the form
$y^2 = x^3 + e$ for some $e\in\GF_q$,
and the twist $F$ of $E_2$ can be written
$cy^2 = dx^3 + e$ for some $c\in\GF_q$ that is not a square
and some $d\in\GF_q$ that is not a cube.
Then the geometric isomorphism $\chi$ can be
taken to be $(x,y) \mapsto (gx,fy)$, where $f,g\in\GFbar_q$
satisfy $g^3 = d$ and $f^2 = c$.
But then it is clear that $\chi$ will
not induce a Galois-equivariant isomorphism $F[\ell]\to E_2[\ell]$.
Therefore, we can tie
$E_1$ and $E_2$ together along their $\ell$-torsion.

This completes the proof of Theorem~\ref{T:distinct1}.
\qed

\section{Proof of Theorem \ref{T:distinct2}}
\label{S:distinct2}

Let $d = t-s$, and suppose $d$ is squarefree.  Then the proof of Corollary~12
of~\cite[p.~1689]{hola} shows that there is no Jacobian isogenous to
$E_1\times E_2$.

On the other hand, suppose there is a prime $\ell$ whose square divides $d$.
The Frobenius $\pi_1$ on $E_1$ is equal to an integer $r$ with $r^2 = q$.
Let $\pi_2$ be the Frobenius on $E_2$, and let $z$ be the element
$(\pi_2 - r)/\ell$ of $\End(E_2)\otimes\Q$.  Then we have
$$z^2 - (d/\ell)z - rd/\ell^2 = 0,$$
so $z$ is integral.  If we replace $E_2$ with an isogenous curve whose
endomorphism ring is maximal, then $z\in\End(E_2)$ so that $\pi_2$ acts as
$r$ on $E_2[\ell]$.  Therefore there are Galois-equivariant anti-isometries
from $E_1[\ell]$ to $E_2[\ell]$.  All of them give rise to Jacobians isogenous
to $E_1\times E_2$, because $E_1$ and $E_2$ are geometrically non-isogenous.
\qed

\section{Proof of Theorem \ref{T:squares}}
\label{S:squares}

Let $\Delta = t^2 - 4q$ and let $R$ be the quadratic order of discriminant
$\Delta$.  Using Serre's appendix to~\cite{lase} or the main results
of~\cite{howe1}, we see that if $\Delta$ is a fundamental discriminant,
then there is a bijection between the set of Jacobians isogenous to
$E\times E$ and the set of indecomposable unimodular hermitian lattices
of rank $2$ over~$R$.  Hoffmann~\cite{hof} shows that if $\Delta$ is
$-3$, $-4$, or $-7$, then there are no such indecomposable unimodular
hermitian lattices, so for these values of $\Delta$ there are no Jacobians
isogenous to $E\times E$.

Suppose $\Delta$ is neither $-3$ nor $-4$ nor $-7$, and suppose $q$ is odd.
Then there is an $E'$ isogenous to $E$ whose automorphism group has order~$2$,
and Corollary~6 of~\cite{hlp} explicitly constructs a genus-$2$ curve whose
Jacobian is isogenous to $E'\times E'$.

Suppose $\Delta$ is neither $-3$ nor $-4$ nor $-7$, and suppose $q$ is a
power of $2$. Then $\Delta\equiv 1\bmod 8$, so $\Delta$ is not the
discriminant of an imaginary quadratic order of class number one (all such
discriminants other than $-7$ are either even or are $5$ modulo $8$).
Therefore there is an elliptic curve $E'$ that is isogenous to $E$ but that
is geometrically non-isomorphic to $E$.  The group schemes $E[2]$ and $E'[2]$
are both isomorphic to the product of $\Bmu_2$ with its dual, so there is an
anti-isometry $E[2]\to E'[2]$.  Since $E$ and $E'$ are geometrically
non-isogenous, this anti-isometry gives rise via Theorem~\ref{T:Kani} to a
Jacobian isogenous to $E\times E'$.
\qed

\begin{remark}
It is also possible to prove the existence of a Jacobian
isogenous to $E^2$ in the case where $t^2-4q \not\in \{-3,-4,-7\}$ directly
from the results of Hoffmann, Serre, and the first author that we cited above,
but some care must be taken in the case where $\Delta$ is not a fundamental
discriminant.
\end{remark}

\section{Twists of Dieudonn\'e modules of supersingular elliptic curves}
\label{S:Dieudonne}

Let $q$ be an even power of a prime $p$, say $q = p^{2a}$.  Our goal in this
section is to prove the following result.

\begin{proposition}
\label{P:SSp-torsion}
Suppose $p>3$.  If $E$ and $E'$ are supersingular elliptic curves over
$\GF_q$ that are not isogenous to one another over $\GF_q$, then the group
schemes $E[p]$ and $E'[p]$ are not isomorphic to one another over $\GF_q$.
\end{proposition}

We will prove this proposition by showing that the Dieudonn\'e modules of
$E[p]$ and $E'[p]$ are not isomorphic to one another.  For concise background
information on Dieudonn\'e modules and $p$-divisible groups,
see~\cite[Ch.~1]{wa} or~\cite[\S3]{oda}.  The first step in our proof of
Proposition~\ref{P:SSp-torsion} will be to compute the twists of the
Dieudonn\'e module of a particular supersingular curve.

Let $E$ be a supersingular elliptic curve over $\GF_q$ whose Weil polynomial
is $(x - \sqrt{q})^2$. As noted at the end of Section~\ref{S:paperintro},
we use $\Q_q$ to denote the unramified extension of $\Q_p$ with residue field
$\GF_q$ and $\Z_q$ to denote the ring of integers of~$\Q_q$.  Let $\sigma$ be
the automorphism of $\Z_q$ over $\Z_p$ that is the lift of the Frobenius
automorphism of $\GF_q$ over $\GF_p$, and let $\frakA$ be the (non-commutative)
ring $\Z_q[F,V]$, where $F$ and $V$ are indeterminates that satisfy
$$FV = VF = p, \quad F\lambda = \lambda^\sigma F,\text{\quad and\quad}
V\lambda^\sigma = \lambda V.$$

Recall that the Dieudonn\'e module $M$ associated to $E$ (or more precisely,
to the $p$-divisible group of $E$) is a certain left $\frakA$-module.
Waterhouse~\cite[p.~539]{wa} computes that $M$ is a free rank-$2$
$\Z_q$-module with a basis $\{x,y\}$ such that
$$Fx = Vx = y \qquad\text{and}\qquad Fy = Vy = px.$$
In particular, we see that all elliptic curves in the isogeny class
with Weil polynomial $(x - \sqrt{q})^2$ have isomorphic Dieudonn\'e
modules.
Waterhouse also notes that the endomorphism ring of $M$ is
isomorphic to the ring of integers $\calO$ of the unique quaternion
algebra $\BH_p$ over~$\Q_p$, and it follows from Waterhouse's analysis
that $M$ gains no further endomorphisms when the base field is extended
to~$\GFbar_q$.

For every $(p^2-1)$st root of unity $\zeta$ in $\Q_{p^2}$ we define a
Dieudonn\'e module $M_\zeta$ as follows:  Let $\xi\in \Q_q$ be a
$(p^{2a}-1)$st root of unity whose norm to $\Q_{p^2}$ is equal to $\zeta$.
Let $M_\zeta$ be a free rank-$2$ $\Z_q$-module generated by
two elements $w$ and $z$, and let $F$ and $V$ act on $w$ and $z$ via
$$Fw = z , Fz = \xi^{-1} p w
    \qquad\text{and}\qquad
  Vw = \xi^{\sigma^{-1}} z, Vz = pw.$$
One can check that this does give a well-defined $\frakA$-module structure
to $M_\zeta$, and that the isomorphism class of $M_\zeta$ does not depend
on the choice of~$\xi$.

\begin{proposition}
\label{P:DMtwists}
The Dieudonn\'e modules $M_\zeta$ over $\GF_q$ are pairwise nonisomorphic
over $\GF_q$, and they are all twists of $M$ over $\GFbar_q$.  The module
$M_\zeta$ is the twist of $M$ by the automorphism $\zeta$ of~$M$.
If $p>3$, then every $\GFbar_q$-twist of $M$ is isomorphic to one of the
$M_\zeta$.
\end{proposition}

\begin{proof}
Let $t$ be an arbitrary element of $M_\zeta$ such that $Ft \not\in p M_\zeta$,
and write $t = aw + bz$ for some $a,b\in \Z_q$.  We see that $a$ must be a unit
of $\Z_q$.  We compute that
\begin{align*}
Ft &= a^\sigma z + b^\sigma \xi^{-1} p w\\
Vt &= a^{\sigma^{-1}} \xi^{\sigma^{-1}} z + b^{\sigma^{-1}} p w.
\end{align*}
It follows that for every $c\in \Z_q$ with
$Ft \equiv c Vt \bmod p M_\zeta$, we have
$$c \equiv  (a^\sigma / a^{\sigma^{-1}}) \xi^{-\sigma^{-1}}\bmod p \Z_q.$$
Taking norms to $\Q_{p^2}$, we find that
$$N_{\Q_q/\Q_{p^2}}(c) \equiv 1\cdot N_{\Q_q/\Q_{p^2}}(\xi^{-\sigma^{-1}})
\equiv \zeta^{-\sigma} \bmod p.$$
Thus, we can recover $\zeta$ from $M_\zeta$, so the $M_\zeta$ are
pairwise nonisomorphic.

Let $B$ be the ring of integers of the maximal unramified extension of $\Q_q$
and let $\frakB = B[F,V]$, where $F$ and $V$ satisfy the same properties as
before.  The base extensions of the $M_\zeta$ to $\GFbar_q$ are
the $\frakB$-modules $\mybar{M}_\zeta$
generated as $B$-modules by $w$ and $z$ and with
$$Fw = z,\quad Fz = \xi^{-1} p w
\quad\text{and}\quad
Vw = \xi^{\sigma^{-1}} z,\quad Vz = pw.$$
Let $\alpha\in B^*$ satisfy $\alpha^{\sigma^2 - 1} = \xi$.  Then one can check
that the map of $\frakB$-modules that sends $x$ to $\alpha w$ and
$y$ to $\alpha^\sigma z$ gives an isomorphism $\phi$
from  $\mybar{M}$ to $\mybar{M}_\zeta$.

The Frobenius automorphism of $\GFbar_q$ over $\GF_q$ acts on
$\Hom(\mybar{M},\mybar{M}_\zeta)$, and we let $\phi^{(q)}$ denote
the image of $\phi$ under this action.
We see that the automorphism $\phi^{-1}\phi^{(q)}$ of $\mybar{M}$
is the map that sends $x$ to $\alpha^{\sigma^{2a} - 1} x$.  Since
$\xi = \alpha^{\sigma^2-1}$, we have
$$\alpha^{\sigma^{2a} - 1} = \xi^{1+\sigma^2+\cdots+\sigma^{2a-2}}
      = \Norm_{\Q_q/\Q_{p^2}}(\xi) = \zeta.$$
Thus, $M_\zeta$ is the twist of $M$ by $\zeta$.

The general theory of twists~\cite{serreCG} shows that the
$\GFbar_q/\GF_q$-twists of $M$ correspond to the elements of the pointed
cohomology set $H^1(\Gal(\GFbar_q/\GF_q), \calO^*)$.  Since the
Galois group acts trivially on $\calO^*$, this cohomology
set consists of the conjugacy classes of $\calO^*$ whose elements
have finite order.

Fix an embedding of $\Q_{p^2}$ into $\BH_p$.  We know
(see~\cite[Thm.~14.5]{reiner}) that there is an element $s\in \BH_p$ with
the properties that $s^2 = p$ and $\BH_p = \Q_{p^2}(s)$, and such that
$s^{-1} x s = x^\sigma$ for all $x\in \Q_{p^2}$.
Suppose that $p>3$, and suppose that $\eta$ is a root of unity in~$\calO$.
Then $\Q_p(\eta)$ is at most a quadratic extension of $\Q_p$, and since
cyclotomic extension of $\Q_p$ have ramification index at least $p-1$ if
they are ramified at all, it follows that $\Q_p(\eta)$ is an unramified
extension of $\Q_p$.   Thus there is a root of unity $\zeta$ in $\Q_{p^2}$
of the same order as $\eta$, and  the Skolem-Noether theorem tells us
that there is an element $x$ of $\BH_p$ that conjugates $\eta$ to $\zeta$.
Let $t$ be the unique power of $s$ such that $tx\in\calO^*$, and let $y =tx$.
Then $y$ conjugates $\eta$ to either $\zeta$ or $\zeta^\sigma$.  Thus,
the elements of $H^1(\Gal(\GFbar_q/\GF_q), \calO^*)$ are represented
by the conjugacy classes that contain roots of unity in $\Q_{p^2}$.  On the
other hand, it is easy to see that no element of $\calO^*$ conjugates
one root of unity in $\Q_{p^2}^*$ to another.
It follows that when $p >3$
there is a bijection between $H^1(\Gal(\GFbar_q/\GF_q), \calO^*)$
and the roots of unity in $\Q_{p^2}$.  Since the $M_\zeta$ are the twists
of $M$ associated to these roots of unity, we find that every twist
of $M$ is isomorphic to some $M_\zeta$.
\end{proof}

\begin{remark}
It is not hard to show that when $p=2$ there are seven
conjugacy classes of roots of unity in $\calO^*$:
in addition to the six conjugacy classes obtained from
the roots of unity in~$\Q_4$,
there is also a single conjugacy class containing a primitive
fourth root of unity.
Likewise, when $p=3$ there are ten
conjugacy classes of roots of unity in $\calO^*$:
eight classes obtained from
the roots of unity in~$\Q_9$,
one class containing a primitive cube root of unity,
and one class containing a primitive sixth root of unity.
\end{remark}

\begin{proof}[Proof of Proposition~\ref{P:SSp-torsion}]
Let $s = p^a$ be the positive square root of $q$.
There are at most 5 isogeny classes of supersingular elliptic
curves over $\GF_q$.  There are always isogeny classes with
Weil polynomials $(x-s)^2$ and $(x+s)^2$.  If $p\equiv 2\bmod 3$,
then there are isogeny classes with Weil polynomial
$x^2 + sx + q$ and $x^2 - sx + q$; each of these isogeny classes
contains two elliptic curves, and they are both twists of the
elliptic curve $y^2 = x^3 - 1$ (by automorphisms of order~$3$
for the former isogeny class, and of order~$6$ for the latter).
If $p\equiv 3 \bmod 4$ then there
is an isogeny class with Weil polynomial $x^2 + q$; there are
two curves in this isogeny classes, each a twist of $y^2 = x^3 - x$
by an automorphism of order~$4$.
(These statements follow from~\cite[Thm.~4.6]{schoof} and its proof.)

We already noted that the Dieudonn\'e module of
every elliptic curve with Weil polynomial $(x-s)^2$ is
isomorphic to the module $M$ defined earlier.  It is also
clear that every elliptic curve with Weil polynomial $(x+s)^2$
has Dieudonn\'e module $M_{-1}$.  When $p\equiv 2 \bmod 3$,
the two curves with Weil polynomial $x^2 + sx + q$ have
Dieudonn\'e modules $M_\zeta$ for two different cube roots of
unity $\zeta$ in the endomorphism ring of $M$, and the
curves with Weil polynomial $x^2 - sx + q$ have
Dieudonn\'e modules $M_\zeta$ for two different sixth roots of
unity.  When $p\equiv 3 \bmod 4$, the two curves with Weil
polynomial $x^2 + q$ have Dieudonn\'e modules isomorphic
to $M_i$ and $M_{-i}$, for a square root $i$ of $-1$ in the
endomorphism ring of~$M$.

Since our elliptic curves $E$ and $E'$ lie in different
isogeny classes, their Dieudonn\'e modules are isomorphic
to $M_\zeta$ and $M_\eta$ for two distinct roots of unity
$\zeta$ and $\eta$ in~$\Z_{p^2}$.  It follows that the
Dieudonn\'e module for $E[p]$
is generated as a $\Z_q$-module by two elements $w$ and $z$ that
satisfy
$$Fw = z,\quad  Fz = 0 \quad\text{and}\quad
Vw = \xi^{1/p} z,\quad  Vz = 0$$
for an element $\xi$ of $\GF_q$ whose norm to $\GF_{p^2}$ is
the reduction of $\zeta$ modulo~$p$.  The same holds for
the Dieudonn\'e module for $E'[p]$, with $\xi$ replaced
by an element $\xi'$ whose norm to $\GF_{p^2}$ is equal to
the reduction of $\eta$ modulo~$p$.

We showed above that $\zeta$ could be recovered from the
module $M_\zeta$.  The same proof shows that $\zeta$ modulo $p$
can be recovered from the Dieudonn\'e module of $E[p]$, and
that $\eta$ modulo $p$
can be recovered from the Dieudonn\'e module of $E'[p]$.
Thus the two Dieudonn\'e modules are not isomorphic to one another,
because the reduction map from roots of unity in $\Z_{p^2}$ to elements
of $\GF_{p^2}$ is injective.
\end{proof}

\begin{remark}
Consider one of the isogeny classes mentioned above
whose Weil polynomial is neither $(x -s)^2$ nor $(x+s)^2$.
It is interesting to note that the two elliptic curves
in this isogeny class have non-isomorphic Dieudonn\'e modules.
It follows that any isogeny between these two
curves must have degree divisible by~$p$.
\end{remark}

\section{Proof of Theorem \ref{T:supersingular}}
\label{S:supersingular}

First suppose that $q$ is not a square.  Because we are assuming that
the characteristic of $\GF_q$ is at least~$5$, there is only one isogeny
class of supersingular elliptic curves over $\GF_q$, and its Weil polynomial
is $x^2 + q$.  From~\cite[Thm~4.5]{schoof} we know that there are
$H(-4p)$ curves in the isogeny class (up to isomorphism over~$\GF_q$),
where $H(\Delta)$ is the Kronecker class number of the discriminant~$\Delta$.
Furthermore,  two curves in the isogeny class are geometrically isomorphic
to one another if and only if they are twists of one another by $-1$,
so the number of distinct $j$-invariants in the isogeny class is $H(-4p)/2$.
In terms of class number of quadratic orders, we have
$$\frac{H(-4p)}{2} = \begin{cases}
                        h(-4p)/2 & \text{if $p\equiv 1\bmod 4$;}\\
                        h(-p)    & \text{if $p\equiv 7\bmod 8$;}\\
                        2h(-p)   & \text{if $p\equiv 3\bmod 8$}.
                     \end{cases}
$$
 From this it follows that when $p\not\in\{5,7,13,37\}$ there are two
curves $E_1, E_2$ in the isogeny class with distinct $j$-invariants.  Since
$E_1[2]$ and $E_2[2]$ are isomorphic Galois modules and $E_1$ and $E_2$
are geometrically non-isomorphic, we can use Theorem~\ref{T:Kani}
to tie $E_1$ and $E_2$ together along their $2$-torsion.

For the remaining cases, we note that if $q$ is an odd power of
a prime $p$ for which $(-2/p) = -1$, the curve $y^2 = x^6 - 5x^4 - 5x^2 + 1$
over $\GF_q$ has Weil polynomial $(x^2 + q)^2$; this is because over $\Q$
its Jacobian is isogenous to the square of the elliptic curve with $j=8000$,
which has complex multiplication by~$\Z[\sqrt{-2}]$.
Since the primes $5$, $7$, $13$, and $37$ all satisfy $(-2/p) = -1$,
we are done.

Now suppose that $q$ is a square, and let $p$ be the unique prime
divisor of $q$.  Recall that there are at most
five isogeny classes of supersingular curves over $\GF_q$; the possible
traces of Frobenius are
\begin{align*}
0             & \text{\qquad if $p\equiv 3 \bmod 4$;}\\
\pm  \sqrt{q} & \text{\qquad if $p\equiv 2 \bmod 3$;}\\
\pm 2\sqrt{q} & \text{\qquad for all $q$.}
\end{align*}

Suppose the traces $s$ and $t$ of our two elliptic curves do \emph{not}
satisfy $s^2 = t^2$.  Then we are to show that there is no Jacobian
isogenous to $E_1\times E_2$.

We begin with a general observation related to Kani's construction
(Theorem~\ref{T:Kani}).   If $E_1$ and $E_2$ are elliptic curves
over $\GF_q$ with traces $s$ and $t$, respectively, and if
$E_1[n]\cong E_2[n]$ as group schemes over $\GF_q$, then
we must have $s \equiv t\bmod n$.  We know that every Jacobian
isogenous to $E_1\times E_2$ is obtained via Kani's construction
for some value of~$n$, and the observation we just made shows
that this value of $n$ must divide~$s-t$.

Suppose that $|s - t| = \sqrt{q}$.  If there were a Jacobian isogenous to
$E_1\times E_2$, it would be attainable through Kani's construction
for some value of $n$ that divides $\sqrt{q}$,
so that this $n$ must be a power of $p$.  But we know from
Section~\ref{S:Dieudonne}
that $E_1[p]\not\cong E_2[p]$, so there are no Jacobians isogenous to
$E_1\times E_2$ in this case.

Suppose that $|s - t| = 2\sqrt{q}$ and that $s\ne -t$,
so that one of $s$ and $t$ is $0$ and the other is $\pm 2\sqrt{q}$.
Say that $s = 0$ and $t = \pm 2\sqrt{q}$.
Note that the endomorphism ring of $E_1$ is isomorphic to
$\Z[i]$, and the Frobenius on $E_1$ is $i\sqrt{q}$;
the Frobenius on $E_2$ is the integer $\pm \sqrt{q}$.
The argument we just gave shows that we cannot obtain a Jacobian
by gluing together $E_1$ and $E_2$ along their $n$-torsion when $n$ is
a multiple of $p$, so if there is a Jacobian isogenous to $E_1\times E_2$
it must be obtained from an anti-isometry $E_1[2]\to E_2[2]$.  But
the Frobenius of $E_2$ acts as a constant on $E_2[2]$, while the
Frobenius of $E_1$ does not, so in particular there are no
anti-isometries from $E_1[2]$ to $E_2[2]$.  Thus there are no
Jacobians isogenous to $E_1\times E_2$.

Suppose that $|s - t| = 3\sqrt{q}$,
so that one of $s$ and $t$ is $\pm\sqrt{q}$ and the other is $\mp 2\sqrt{q}$.
Say that $s = \pm\sqrt{q}$ and $t = \mp 2\sqrt{q}$.
Note that the endomorphism ring of $E_1$ is isomorphic to
$\Z[\omega]$ for some cube root of unity $\omega$,
and the Frobenius on $E_1$ is $\mp\omega\sqrt{q}$;
the Frobenius on $E_2$ is the integer $\mp \sqrt{q}$.
Again we see that
we cannot obtain a Jacobian
by gluing together $E_1$ and $E_2$ along their $n$-torsion when $n$ is
a multiple of $p$, so if there is a Jacobian isogenous to $E_1\times E_2$
it must be obtained from an anti-isometry $E_1[3]\to E_2[3]$.  But
the Frobenius of $E_2$ acts as a constant on $E_2[3]$, while the
Frobenius of $E_1$ does not, so again we see there are no
Jacobians isogenous to $E_1\times E_2$.

Now suppose we have $s^2 = t^2$.  We must show that there is a Jacobian
isogenous to $E_1\times E_2$. There are three cases to consider.

\subsubsection*{The case  $s^2 = t^2 = q$.}
This case arises only
when $ p \equiv 2\bmod 3$.  Let $a$ be a generator of~$\GF_q^*$, and
consider the curve $C$ defined by $y^2 = x^6 + a$.  Arguing as
in~\cite[\S3]{hlp}, we see that the Jacobian of $C$ is isogenous to
$F_1\times F_2$, where $F_1$ is the elliptic curve
$y^2 = x^3 + a$ and $F_2$ is the elliptic curve $y^2 = x^3 + a^2$.
Let $F_0$ be the elliptic curve $y^2 = x^3 + 1$,  let $b$ be a
sixth root of $a$ in $\GFbar_q$, let $\zeta$ be the primitive
sixth root of unity $b^{q-1}$, and let $\omega$ be the
order-$6$ automorphism $(x,y)\mapsto(\zeta^2 x, \zeta^3 y)$
of $E_0$.
Since $F_0$ is defined over $\GF_p$, its Frobenius endomorphism
over $\GF_q$ is either $\sqrt{q}$ or $-\sqrt{q}$.
It is easy to see that $F_1$ is the twist of $F_0$ by $\omega$
and that $F_2$ is the twist of $F_0$ by $\omega^2$, and it follows
that $F_1$ and $F_2$ have traces of opposite sign, and they
both are square roots of~$q$.

Similar reasoning shows that the Jacobian of the curve
$y^2 = x^6 + a^2$ is isogenous to either $F_1\times F_1$
or $F_2\times F_2$; furthermore, whichever product of
elliptic curves we get from $y^2 = x^6 + a^2$, we get the
other product from the quadratic twist
$ay^2 = x^6 + a^2$.

Thus, whenever $s^2 = t^2 = q$ there is a Jacobian
isogenous to $E_1\times E_2$.

\subsubsection*{The case  $s^2 = t^2 = 0$.}
This case occurs only when $p \equiv 3\bmod 4$.  Let $F_0$ be the elliptic
curve over $\GF_q$ defined by $y^2 = x^3 - x$, so that $j(F_0) = 1728$ and
$F_0$ has an automorphism $i$ of order $4$.  The two elliptic curves $F_1$
and $F_2$ over $\GF_q$ with trace $0$ are the twists of $F_0$ by $i$ and
by~$-i$.  Let us fix, once and for all, two $\GFbar_q$-isomorphisms
$\phi_1\colon F_0\to F_1$ and $\phi_2\colon F_0\to F_2$.  Using these
isomorphisms, we will identify geometric points of $F_1$ and $F_2$ with
geometric points of $F_0$, and we will identify the geometric automorphism
groups of $F_1$ and $F_2$ with $\Aut F_0$.  Let $s$ be the positive square
root of $q$, and reindex the curves if necessary so that the $q$-power
Frobenius on $F_1$ is equal to~$si$.  Then the Frobenius on $F_2$ is equal
to~$-si$.

We will show that there is a Jacobian isogenous to $F_1\times F_2$
by gluing the two curves together along their $4$-torsion subgroups,
as in Theorem~\ref{T:Kani4}.

Let $P$  be a geometric point of $F_0$ such that $P$ and $iP$
generate $F_0[4]$.  Let $Q = iP$.  Let $\psi$ be the isomorphism
$F_1[4]\to F_2[4]$ that sends $\phi_1(P)$ to $\phi_2(P+2Q)$
and $\phi_1(Q)$ to $\phi_2(2P-Q)$.  It is easy to check that
$\psi$ is a Galois-equivariant anti-isometry with respect to
the Weil pairing.  We will be finished if we can show that
neither condition (a) nor condition (b) of Theorem~\ref{T:Kani4}
holds.

We know that $F_1$ and $F_2$ are both elliptic curves with
$j$-invariant~$1728$.  Let $\Phi_3(j,j')\in\Z[j,j']$ be
the classical modular polynomial for $3$-isogenies.
We compute that
$$\Phi_3(1728,1728) = 2^{36} \cdot 3^6 \cdot 7^8\cdot 11^4,$$
so when $p>11$ there are no geometric $3$-isogenies from $F_1$ to $F_2$,
and condition (a) of Theorem~\ref{T:Kani4} does not hold.
For $p=7$ and $p=11$, we can explicitly write down all $3$-isogenies
from $F_0$ to $F_0$ and note that none of them induce the given
anti-isometry $\psi$ from $F_1[4]$ to $F_2[4]$.

If condition (b) of Theorem~\ref{T:Kani4} were to hold,
there would be two order-$2$ subgroups $G_1$ and $G_2$ of $F_0[4]$
and an automorphism $\alpha$ of $F_0$ such that
\begin{align*}
P + \alpha(P+2Q) \in G_1\quad &\quad P - \alpha(P + 2Q) \in G_2\\
Q + \alpha(2P-Q) \in G_1\quad &\quad Q - \alpha(2P - Q) \in G_2
\end{align*}
Note that the only automorphisms of $F_0$ are $\pm1$ and $\pm i$.
We check that $P+\alpha(P)$ is a $2$-torsion element only when
$\alpha = \pm 1$.  But if $\alpha=1$ then we find that $G_1$
contains both $2P+2Q$ and $2P$, a contradiction, while if
$\alpha=-1$ then $G_1$ contains both $2Q$ and $2P+2Q$, another
contradiction.  Thus, condition (b) of Theorem~\ref{T:Kani4} does not hold.
It follows that there is a Jacobian with Weil polynomial
$(x^2 + q)^2$.

\subsubsection*{The case  $s^2 = t^2 = 2q$.}
Note that the Galois group  of $\GFbar_q/\GF_q$ acts trivially
on both $E_1[2]$ and $E_2[2]$, because for each curve
the Frobenius endomorphism is either $\sqrt{q}$ or $-\sqrt{q}$.
There are therefore six Galois-equivariant anti-isometries
$E_1[2]\to E_2[2]$.  Since we are not in characteristic $2$ or~$3$,
the number of isomorphisms $E_1\to E_2$ is at most~$6$, and in the
case that there are $6$ isomorphisms, there are only $3$ induced
isomorphisms $E_1[2]\to E_2[2]$.  Thus, at least one of the
six anti-isometries $E_1[2]\to E_2[2]$ is not reducible, so there
is a Jacobian isogenous to $E_1\times E_2$. \qed


\part{Simple supersingular abelian surfaces as Jacobians}
\label{Part:simple}

\section{Introduction}
\label{S:IntroSimple}

Let $k=\fq$ be a finite field of characteristic $p>3$, and let $\calA$
be an isogeny class of simple supersingular abelian surfaces over $k$ that
split over the quadratic extension of~$k$.  In this part of the paper
we determine whether or not there is a Jacobian in~$\calA$.

In \cite[Table 1]{mn} we find a list of all simple supersingular isogeny
classes of abelian surfaces defined over a finite field, with an indication
of the smallest field extension over which each class splits.  We present in
Table~\ref{Table:simpleSS} the isogeny classes over finite fields $\fq$ of
characteristic $p>3$ that split over~$\ff2$.  Theorem~\ref{T:PPAS} shows
that these classes are all principally polarizable.  We will show that all
of these isogeny classes contain Jacobians, except for one special case.

\begin{table}
\begin{center}
\as{1.1}
\begin{tabular}{|l|l|l|}
\hline
 $(a,b)$  & Conditions on $p$        & Conditions on $q$ \\
\hline\hline
$(0,0)$   & $p \not\equiv 1 \bmod 4$ & ---               \\ \hline
$(0,q)$   & $p \not\equiv 1 \bmod 3$ & $q$ nonsquare     \\ \hline
$(0,-q)$  & $p \not\equiv 1 \bmod 3$ & ---               \\ \hline
$(0,-2q)$ & ---                      & $q$ nonsquare     \\ \hline
$(0,2q)$  & $p \equiv 1 \bmod 4$     & $q$ square        \\ \hline
\end{tabular}
\end{center}
\vspace{1ex}
\caption{Simple supersingular isogeny classes over $\fq$
that split over $\ff2$.  The characteristic $p$ of $\fq$ is
assumed to be at least~$5$.}
\label{Table:simpleSS}
\end{table}

\begin{theorem}
\label{T:part2}
Let $\calA_{(a,b)}$ be an isogeny class of simple supersingular abelian surfaces
over a finite field $\fq$ of characteristic $p>3$.  Then
$\calA_{(a,b)}$ does not contain a Jacobian if and only
if $q$ is a square, $p\equiv 11\bmod 12$, and $(a,b) = (0,-q)$.
\end{theorem}

This part of the paper is organized as follows.  In Section~\ref{one} we review
results of Oort~\cite{oort}, Katsura and Oort~\cite{oort2}, and Ibukiyama,
Katsura, and Oort~\cite{IKO} on supersingular abelian surfaces over the
algebraic closure of a finite field, paying special attention to the principal
polarizations of these surfaces. In Section~\ref{two} we look at supersingular
surfaces over finite fields and determine which of their geometric principal
polarizations can be defined over the base field.
Finally, in Section~\ref{three} we use the results of Section~\ref{one}
and~\ref{two}, together with some explicit constructions, to prove
Theorem~\ref{T:part2}.

\section{Supersingular surfaces, quaternion lattices, and polarizations}
\label{one}

In this section we review some results of Oort~\cite{oort},
Katsura and Oort~\cite{oort2}, and Ibukiyama, Katsura, and Oort~\cite{IKO}
on supersingular abelian surfaces, quaternion hermitian forms,
and polarizations.  The results on abelian surfaces assume that the base field
is algebraically closed; we will consider the case of finite base fields
in Section~\ref{two}.

\subsection{Supersingular abelian surfaces}
\label{SS:sas}

Let $E$ be an elliptic curve over $\fp$ with trace~$0$, so that $E$ is
supersingular and all of the geometric endomorphisms of $E$ are
defined over~$\fpt$.  Let $K$ be the algebraic closure of $\fp$ and
let $\calO$ be the $K$-endomorphism ring of $E$; the algebra $B=\calO\otimes\Q$
 is a definite quaternion algebra over~$\Q$ with discriminant~$p$,
and $\calO$ is a maximal order in~$B$.
We will denote the canonical anti-involutions of $\calO$ and $B$ by
$x\mapsto\xbar$.

Let $\pi$ denote the $p$-power Frobenius endomorphism on $E$,
and fix a $K$-isomorphism between $E[\pi]$ and $\alpha_p$, where
$\alpha_p$ is the unique local-local group scheme over $\fp$;
then we can identify $\Hom_K(\alpha_p,E)$ with $\End_K(\alpha_p)=K.$
The kernel of the restriction map
\begin{align*}
\widetilde{\ }:\End_K(E)&\to     \End_K(\alpha_p)\\
                     u &\mapsto \tilde{u}=u|_{\alpha_p}
\end{align*}
is a two-sided prime ideal $\frakP$ of $\calO$ above $p$, with residual
degree~$2$. The restriction map thus gives a natural embedding
$\calO/\frakP\hookrightarrow \End_K(\alpha_p)=K$ with image $\fpt$.
Since $\pi^2=-p$, the prime ideal $\frakP$ is principal and generated by $\pi$.

For every $(i,j)\in K^2$, we denote by $A_{ij}$ the abelian surface over $K$
given by the following diagram:
$$0 \to \alpha_p \longmapright{(i,j)} E\times E \to A_{ij} \to 0.$$
It is easy to check that
\begin{align*}
\label{aij}
A_{ij}=A_{i'j'}&\sii (i,j)(\alpha_p)=(i',j')(\alpha_p)\\
               &\sii \exists a\in K^* \text{\ such that\ }(i',j')=a(i,j).
\end{align*}
Thus, the set of all $A_{ij}$ (apart from $A_{00}=\et$) is parameterized
by~$\pr1(K)$.

For every $i\in K$ the composition
$\alpha_p \mapright{i} E \mapright{u} E$ corresponds to the
element $i\,\tilde{u}$ of~$K$.
For every endomorphism $\alpha\in\End_K(\et)\cong M_2(\calO)$ and every
$[i:j]\in\pr1(K)$, the composition
$$\alpha_p \longmapright{(i,j)} \et \mapright{\alpha}\et$$
has the same image as the element $\widetilde{\alpha}[i:j]\in\pr1(K)$, where
$\widetilde{\alpha}\in M_2(\fpt)$ is obtained by reduction modulo $\frakP$ of
the entries of $\alpha$, and the action of $M_2(\fpt)$ on $\pr1(K)$ is the
usual projective action.

If $A$ is an abelian surface over $K$ we denote by $a(A)$ the
quantity
$$a(A)=\dim \Hom_K(\alpha_p,A),$$
sometimes called the \emph{$a$-number} of $A$.
When $A$ is a supersingular abelian surface we have $a(A)\in\{1,2\}$.
The value of $a(A)$ gives us information about the global structure
of $A$, as the following result shows.

\begin{proposition}
\label{P:avalue}
We have $a(A)=2$ if and only if $A\cong E\times E$.
We have $a(A) = 1$ if and only if $A\cong A_{ij}$
for some $[i:j] \in \pr1(K)\setminus\pr1(\fpt)$.
Furthermore, if $a(A) = 1$ then $a(A/\alpha_p)=2$.
\end{proposition}

\begin{proof}
This follows from [Oort 75, Introduction], [Oort 75, Thm.~2],
and [Oort 75, Cor.~7].
\end{proof}

\subsection{Quaternion hermitian forms and lattices}

Most of the material that we present without reference in this section
can be found in~\cite{shimura}.

Let $B$ be a definite quaternion algebra over~$\Q$ with discriminant $p$.
There is a positive definite hermitian form on the right $B$-module $B^2$,
which is unique up to base change over $B$; it is given explicitly
by $\sum \mybar{x_i} y_i$, where
$x\mapsto\xbar$ is the standard involution on~$B$.
For every prime $\ell$ (possibly equal to $p$) we set
$B_{\ell}=B \otimes \Q_{\ell}$.  Then the hermitian form on $B^2$
extends to give a hermitian form on $B_\ell^2$.
Let $\dag$ denote the conjugate-transpose involution on $M_2(B)$
and on $M_2(B_\ell)$, where `conjugation' means the standard involution.
Then the groups of similitudes of the hermitian forms on $B^2$ and
on $B_\ell^2$ are given by
\begin{align*}
G    &=\{g \in M_2(B) \tq
         g^{\dag}g=n(g) I \text{\ for some\ }
         n(g) \in \Q^*\}\\
\intertext{and}
G_{\ell}&=\{g \in M_2(B_{\ell}) \tq
         g^{\dag}g=n(g) I\text{\ for some\ }
         n(g) \in \Q_{\ell}^*\}.
\end{align*}

Let $\calO$ be a maximal order of $B$. A $\Z$-lattice $L$ in
$B^2$ is called a \emph{{\rm(}right\/{\rm)} $\calO$-lattice} if $L$ is a right
$\calO$-module. Two $\calO$-lattices $L_1$ and $L_2$ are
\emph{globally equivalent} if $L_1=g L_2$ for some $g\in G$,
and are \emph{locally equivalent at $\ell$}
if $L_1 \otimes \Z_\ell=g (L_2 \otimes \Z_\ell)$ for some $g \in G_\ell$.
One denotes by $\Aut(L)=\{g \in G \tq g L=L\}$ the
\emph{automorphism group} of $L$ and by $\Aut'(L)=\Aut(L)/{\pm 1}$
the \emph{reduced automorphism group} of $L$.  These groups are finite
because the hermitian form on $B^2$ is positive definite.

A \emph{genus} of $\calO$-lattices is a set of $\calO$-lattices in $B^2$ that
are equivalent to one another locally at every prime $\ell$.  There are only
two genera: the \emph{principal genus} $\calL_2(p,1)$ that contains the right
$\calO$-lattices in $B^2$ that are equivalent to $\calO_\ell^2$ for all $\ell$,
and the \emph{non-principal genus} $\calL_2(1,p)$ that contains the right
$\calO$-lattices in $B^2$ that are equivalent to $\calO_\ell^2$ for all
$\ell \ne p$ and equivalent at $p$ to
\begin{equation}
\label{elp}
\xi \left( \begin{array}{cc} 1 &0 \\ 0 & \pi  \end{array} \right)
\calO_p^2, \qquad \mbox{ where $\xi\in\gl2{\calO_p}$ satisfies }
\xi^{\dag}\xi=\left( \begin{array}{cc} 0 &1 \\ 1 & 0  
\end{array} \right),
\end{equation}
with  $\pi$  a prime element of  $\calO_p$.
One denotes by $H_2(p,1)$ the (finite) number of global equivalence classes
in $\calL_2(p,1)$, and by $H_2(1,p)$ the number of global equivalence
classes in $\calL_2(1,p)$.

On the other hand one can define two special sets of positive definite
hermitian matrices. Let $\frakP$ be the two sided prime ideal of
$\calO$ above $p$.

\begin{definition}
We define $\Lambda^{\prc}$ to be the set of matrices
$H$ in $\gl2{\calO}$ such that
$$H=\left( \begin{matrix}s & r \\ \mybar{r} & t \end{matrix}\right)
\text{ with } st-r \mybar{r}=1,$$
where $s$ and $t$ are positive integers.
We define $\Lambda^{\nprc}$ to be the set of matrices
$H$ in $M_2(\calO)$ such that
$$H=\left( \begin{matrix} p s & r \\ \mybar{r} & p t \end{matrix}\right)
\text{ with } p^2st-r \mybar{r}=p,$$
where $s$ and $t$ are positive integers and where $r\in\frakP$.
\end{definition}

Two matrices $H_1,H_2$ that both lie in $ \Lambda^{\prc}$ or in
$\Lambda^{\nprc}$ are said to be \emph{equivalent} if there exists an
$\alpha \in \GL_2(\calO)$ such that  $\alpha^{\dag} H_1 \alpha=H_2$.
For $H$ in $ \Lambda^{\prc}$ or $\Lambda^{\nprc}$, we let
$$\Aut(H)=\{\alpha \in \GL_2(\calO) \mid \alpha^{\dag} H \alpha=H\}$$
be the \emph{automorphism group} of $H$ and $\Aut'(H)=\Aut(H)/\pm 1$ the
\emph{reduced automorphism group} of $H$. These groups are again finite.

One can relate lattices and hermitian forms in the following way.

\begin{proposition}
\label{corr}
There are bijective correspondences
\begin{align*}
\left\{\vcenter{\hsize=4.5cm \noindent
global equivalence classes\hfil\break of lattices in $\calL_2(p,1)$}\right\}
&\longleftrightarrow
\left\{\vcenter{\hsize=3.5cm \noindent
equivalence classes\hfil\break of matrices in $\Lambda^{\prc}$}\right\}\\
\intertext{and}
\left\{\vcenter{\hsize=4.5cm \noindent
global equivalence classes\hfil\break of lattices in $\calL_2(1,p)$}\right\}
&\longleftrightarrow
\left\{\vcenter{\hsize=3.5cm \noindent
equivalence classes\hfil\break of matrices in $\Lambda^{\nprc}$}\right\}
\end{align*}
that preserve automorphism groups and reduced automorphism groups,
as abstract groups.
\end{proposition}

\begin{proof}
The bijective correspondences are provided by Lemmas~2.3 and~2.5
of~\cite{IKO} in the principal case, and by Lemmas~2.6 and~2.7
of~\cite{IKO} in the non-principal case.  The fact that the
bijections provided by these lemmas preserve automorphism groups is
easily seen from the proofs of the lemmas.
\end{proof}

The automorphisms $\alpha$ of the hermitian forms in $\Lambda^{\nprc}$ are
determined (up to $\pm1$) by the projective action of $\widetilde{\alpha}$.
This fact is probably well-known, but for lack of a suitable reference we
include a short proof.

\begin{lemma} \label{reduc}
Let $\Gamma\subseteq \gl2{\calO}$ be a finite subgroup of order prime to~$p$.
Then reduction modulo $\frakP$ determines an embedding
$\widetilde{\phantom{m}}\colon\Gamma\hookrightarrow \gl2{\fpt}$.
\end{lemma}

\begin{proof}
If $\alpha\in\Gamma$ is the identity modulo $\frakP$, then
$\alpha$ is an element of the multiplicative group $1+\pi M_2(\calO)$.
This group has no torsion element of order prime to~$p$, so we must have
$\alpha=1$.

More explicitly, every element of the group $1+\pi M_2(\calO)$ can be written in
the form  $1+\pi^NM$, with $N>0$ and $M\not\in M_2(\frakP)$. Thus, for every
positive integer $n$ we have
$$
(1+\pi^NM)^n\equiv 1+n\pi^NM \bmod \mathfrak{P}^{2N},
$$
and if $n$ is prime to $p$ we cannot have $(1+\pi^NM)^n=1$, since
$n\pi^NM \not\equiv0\bmod \mathfrak{P}^{2N}$.
\end{proof}

\begin{proposition}\label{det1}
Let $H$ be an element of $\Lambda^{\nprc}$ and let $\Gamma\subseteq\Aut(H)$ be
a subgroup of order prime to $p$. Then reduction modulo $\mathfrak{P}$
gives embeddings
$$\Gamma\hookrightarrow\SL_2(\fpt)\text{\quad and\quad}
  \Gamma/\{\pm1\}\hookrightarrow\pgl2{\fpt}.$$
\end{proposition}

\begin{proof}
By the above lemma, we need only check that $\det(\widetilde{\alpha})=1$
for all $\alpha\in\Aut(H)$.  Put $H=x^{\dag}x$ for some $x\in\gl2{B}$ and
let $L=x\calO^2$ be the lattice in the non-principal genus attached to $H$
as in Proposition~\ref{corr}. Since the lattice $L\otimes\Z_p$ is equivalent
to the lattice given in~\eqref{elp}, there exist $g\in G_p$ and
$\beta\in \gl2{\calO_p}$ such that
$x=g\xi\begin{pmatrix}1&0\\0&\pi\end{pmatrix}\beta$. We compute that
$$
H = x^{\dag}x
  = \beta^{\dag}
    \begin{pmatrix}1&0\\0&-\pi\end{pmatrix}
    \xi^{\dag}g^{\dag}g\xi
    \begin{pmatrix}1&0\\0&\pi\end{pmatrix}
    \beta
  = n(g)\beta^{\dag}
    \begin{pmatrix}0&\pi\\-\pi&0\end{pmatrix}
    \beta.
$$
Let $\gamma=\beta\alpha\beta^{-1}$.  Then from $H=\alpha^{\dag}H\alpha$
we find that
\begin{equation}\label{gap}
\gamma^{\dag}
  \begin{pmatrix}0&\pi\\-\pi&0\end{pmatrix}
  \gamma
 =\begin{pmatrix}0&\pi\\-\pi&0\end{pmatrix}.
\end{equation}
Since $\gamma$ and $\beta$ lie in $\gl2{\calO_p}$, we can reduce
them modulo $\frakP$ as well; thus, it is sufficient to check that
$\det(\widetilde{\gamma})=1$. Now, $\overline{\gamma}\pi=\pi\gamma'$
with $\widetilde{\gamma'}=\widetilde{\gamma}$. Hence, we can cancel
$\pi$ in both sides of \eqref{gap} to get
$$
(\gamma')^t\begin{pmatrix}0&1\\-1&0\end{pmatrix}
  \gamma
  =\begin{pmatrix}0&1\\-1&0\end{pmatrix},
$$
and this implies that $\widetilde{\gamma}$ belongs to the symplectic group
and has determinant equal to~$1$.
\end{proof}

Katsura and Oort determined the groups that can occur as the reduced
automorphism group of a hermitian form in $\Lambda^{\nprc}$ (see
\cite{oort2} and \cite[Lem.~2.1]{ibu}). This result will play a
crucial role in our strategy.

\begin{theorem}
\label{mass1}
If $p\ge7$, then the reduced automorphism group of a
hermitian matrix in  $\Lambda^{\nprc}$ is isomorphic as an abstract group to
one of the following groups{\/\rm:}
$$
\Z/n\Z \text{\ for some $n\in\{1,2,3\}$};
\quad D_{2n} \text{\ for some $n\in\{2,3,6\}$};
\quad A_4;
\quad S_4;
\quad A_5.
$$
If $p=3$ or $p=5$, then $H(1,p)=1$ and the reduced
automorphism group of the single class is isomorphic
to $A_6$ when $p=3$ and to $\pgl2{\mathbb{F}_5}$ when $p=5$.
\end{theorem}

Given a monic polynomial $f\in\Z[x]$ and a subgroup $\Gamma$ of
$\GL_2(B)$, let $\Gamma_f$ denote the set of elements of $\Gamma$
whose reduced characteristic polynomials (as elements of $M_2(B)$)
are equal to $f$.  For each possible reduced automorphism group
$\Gamma'$ of a hermitian matrix in $\Lambda^{\nprc}$,
Ibukiyama~\cite[Thm.~7.1]{ibu} determined the cardinality of the
set of equivalence classes of hermitian matrices
$H\in \Lambda^{\nprc}$ with $\Aut'(H)\cong \Gamma'$.
An important ingredient in this computation is the
determination of mass formulas for the number of elements in $\Gamma_f$
for all $\Gamma$ and~$f$.
Given a monic degree-$4$ polynomial $f\in\Z[x]$,
we define $m(f)$ to be the quantity
$$m(f):=\sum_{i=1}^h\frac{|\Gamma_{i,f}|}{|\Gamma_i|},$$
where $h=H(1,p)$ is the number of classes in $\Lambda^{\nprc}$
and where the $\Gamma_i$ are the automorphism groups of a set
of representatives for the equivalence classes of~$\Lambda^{\nprc}$.
Ibukiyama computed the masses $m(f)$ explicitly
(\cite[Thm.~2.2]{ibu}).

\begin{theorem}
\label{mass2}
Assume that $p\ge 7$. Then $m(f)=0$ for all polynomials $f\in\Z[x]$ except for
those with $f(x)$ or $f(-x)$ belonging to the following list{\/\rm:}
$$\as{1.1}
\begin{array}{lll}
f_1=(x-1)^4, &\qquad f_2=(x^2+1)^2,       &\qquad f_3=(x^2+x+1)^2,\\
f_4=x^4+1,   &\qquad f_5=x^4+x^3+x^2+x+1, &\qquad f_6=x^4-x^2+1.
\end{array}
\as{1}
$$
Moreover, $m(f_2)>0$ for all $p\ge 7$, and
\begin{align*}
m(f_4)>0 \ &\mbox{ if and only if }\ p\equiv 3,\,5\bmod8,\\
m(f_6)>0 \ &\mbox{ if and only if }\ p\equiv 5\bmod{12}.
\end{align*}
\end{theorem}

\begin{remark}
There is a similar result in the principal genus case;
see~\cite[Part~I]{hi1}.
\end{remark}

\subsection{Polarizations}
\label{pohe}
Later in the paper we will need to understand the principal polarizations
on the supersingular abelian surfaces over the algebraic closure $K$ of a finite
field $\fq$.  In this section we present the relevant results.

Recall that in Section~\ref{SS:sas} we chose a trace-$0$
elliptic curve $E$ over~$\fp$.
The $K$-endomorphism ring $\calO$ of $E$ is a maximal order in the
quaternion algebra $B$ over $\Q$ with discriminant $p$.

Let $\lambda_0$ be the product principal polarization on
$\et$ and let $\dag$ be the Rosati involution on
$\End_K(\et)$ associated to this polarization.
It is well-known that under the natural isomorphism
$\End_K(\et) \cong M_2(\calO)$, the Rosati
involution becomes the conjugate-transpose involution.

The polarization $\lambda_0$ induces an injection from the N\'eron-Severi
group $\NS(\et)$ to $\End_K(\et)$ by
$\lambda \mapsto \lambda_0^{-1}\lambda$.  The image of this map is also
well-known (see for instance \cite[Prop.~2.8]{IKO}).

\begin{proposition}
The  map given above induces a bijection between $\NS(\et)$ and the set
of hermitian matrices in $M_2(\calO)$.  Moreover, this map restricts to a
bijection between the set of principal polarizations on $\et$ and
$\Lambda^{\prc}$.
\end{proposition}

We can understand in a similar way the principal polarizations of the
supersingular surfaces that are not geometrically isomorphic to~$\et$.
Let $A$ be a supersingular abelian surface with $a(A) = 1$, and let
$\psi$ be the natural degree-$p$ isogeny from $\et$ to $A$
(see \S\ref{SS:sas}).
Then we can define a map $\NS(A) \to M_2(\calO)$ by
$$\lambda \mapsto \lambda_0^{-1} \hat{\psi} \lambda \psi.$$
Proposition~2.14 of~\cite{IKO} tells us the following.

\begin{proposition}\label{P:nonprpol}
The map given above induces a bijection between the set of principal
polarizations on $A$ and the set $\Lambda^{\nprc}$.
\end{proposition}

\section{Supersingular surfaces over finite fields.}
\label{two}
Let $k=\fq$ be a finite field of characteristic $p$ that has even
degree over $\fp$, and let $K=\GFbar_q$.
In this section we answer some basic questions concerning
supersingular abelian surfaces over~$k$, their isogeny classes, and their
principal polarizations.

Suppose $A$ is a supersingular abelian surface over $k$.  If $a(A) = 2$
then $A$ is a $K/k$-twist of the abelian surface~$\et$, where,
as before, $E$ is a trace-$0$ elliptic curve over $\fp$.  On the other hand,
if $a(A) = 1$ then there is a unique copy of $\alpha_p$ in $A$, which must
necessarily be defined over $k$, and by Proposition~\ref{P:avalue}
the quotient $A/\alpha_p$ has $a$-number $2$.  Therefore, every $A$
is either a $K/k$ twist of $\et$ or a quotient of such a
 twist by a rank-$p$ subgroup. In particular,
every isogeny class of supersingular surfaces over $k$ contains
a $K/k$-twist of $\et$.

Thus, to understand the supersingular abelian surfaces over $k$ and
their principal polarizations,
we need only answer the following questions:
\begin{itemize}
\item What are the $K/k$-twists of $\et$,
      and what are the Weil polynomials of these twists?
\item Which rank-$p$ geometric subgroups of these twists
      can be defined over $k$?
\item Which geometric polarizations of these twists
      can be defined over $k$?
\end{itemize}
In this section we will answer these questions.

\begin{remark}
We could ask the same
questions for arbitrary finite fields instead of limiting ourselves
to those that contain $\fpt$, but we will only need the answers for
even-degree extensions of $\fp$, and the answers for odd-degree extensions
of $\fp$ are slightly more awkward to state.  The answers are simpler
for the fields that contain $\fpt$ because for these fields the Galois
group of $K/k$ acts trivially on $\End_K(\et)$.
\end{remark}

The first of our three questions is easy to answer.  We know that the twists
of $\et$ correspond to elements of the cohomology set
$H^1(\Gal(K/k),\Aut_K(\et))$, and since all of the
geometric endomorphisms of $E$ are defined over $k$,
this cohomology set consists of the conjugacy classes of
the elements of finite order in $\Aut(\et)$. If $A$ is a twist of $\et$,
and if $f\colon\et\to A$ is a geometric isomorphism, then $\alpha:=f^{-1} f^\sigma$
is the automorphism of $\et$ that corresponds to
the twist~$A$; here $\sigma$ is the Frobenius automorphism of $K/k$.

Let $\pi$ be the $q$-power Frobenius on $\et$ and let $\pi_A$ be the
$q$-power Frobenius of a twist $A$ of $\et$. The pullback of $\pi_A$
via the geometric isomorphism $f$ is equal to $\alpha\pi$, so that
$\pi_A$ and $\alpha\pi$ have the same characteristic polynomial.
Since $\pi=\pm\sqrt{q}$ is an integer in $\End_K(\et)$, the
characteristic polynomial of $\pi_A$ is $\pi^4 h(x/\pi)$, where
$h\in\Z[x]$ is the characteristic polynomial of $\alpha$.
The same argument is valid in a more general situation:

\begin{proposition}
\label{P:charpoly}
Let $A$ and $B$ be abelian surfaces over $k$ with $q$-power Frobenius
endomorphisms $\pi_A$ and $\pi_B$, respectively.
Let $f\colon B\longrightarrow A$ be a $K$-isomorphism and let $h\in\Z[x]$
be the characteristic polynomial of $\alpha=f^{-1} f^\sigma\in\Aut_K(B)$.
If $\pi_B$ acts as an integer on $B$ then the characteristic polynomial
of $\pi_A$ is $\pi_B^4 h(x/\pi_B)$.
\end{proposition}

We turn to the second question.
Given $(i,j)\in K^2$ with $(i,j)\neq(0,0)$, we would like to know whether the
subgroup $f((i,j)(\alpha_p))$ of $A$ is definable over $k$.

\begin{proposition}
\label{P:subgroups}
The subgroup $f((i,j)(\alpha_p))$ of $A$ is definable over $k$
if and only if $[i:j]$ and $\widetilde{\alpha}[i^\sigma:j^\sigma]$
are equal in $\pr1(K)$.
\end{proposition}

\begin{proof}
The morphism $f\circ(i,j)\colon \alpha_p\longrightarrow A$ is defined
over $k$ if and only if it is invariant under the action of the Galois
group of $K/k$, so we have
$$
f\circ(i,j)=f^{\sigma}\circ(i^{\sigma},j^{\sigma})
   \sii  (i,j)=\alpha\circ(i^{\sigma},j^{\sigma})
   \sii [i:j]=\widetilde{\alpha}\,[i^{\sigma}:j^{\sigma}].
$$
\end{proof}

It will also be useful to know when we can be assured of
the existence of a rational local-local subgroup of $A$
that gives rise to a quotient with $a$-number $1$.
By Proposition~\ref{P:avalue} we
need $[i:j]\in\pr1(K)\setminus\pr1(\fpt)$ such that the morphism
$f\circ(i,j)$ is defined over $k$.

\begin{proposition}
\label{P:goodsubgroups}
Let $H\in\Lambda^{\nprc}$ and let $\alpha\in \Aut(H)$ be
an automorphism of order not divisible by $p$.  Then there exists
an element $[i:j]\in\pr1(K)\setminus\pr1(\fpt)$
such that $[i:j] = \widetilde{\alpha}[i^\sigma:j^\sigma]$,
unless $q = p^2$ and
$\alpha = \pm1$.
\end{proposition}

\begin{proof}
The equation $[i:j] = \widetilde{\alpha}[i^q:j^q]$
can be rewritten as a homogeneous equation in $i$ and $j$
of degree $q+1$, and it is easy to verify that
the subscheme of $\pr1$ defined by this equation
is nonsingular, so there are $q+1$ points $[i:j]$
that satisfy the equation.  If $q>p^2$, we are guaranteed
a solution that does not lie in $\pr1(\fpt)$.
If $q=p^2$
and every element of $\pr1(\fpt)$ is a root of the
equation, then we see that
$\widetilde{\alpha}$ fixes every element of $\pr1(\fpt)$, and
Proposition~\ref{det1} shows that $\alpha=\pm1$.
\end{proof}

The third of our three questions asks when a geometric polarization
of $\et$ gives rise to a
polarization of $A$ defined over $k$.

\begin{proposition}
\label{P:polarizations}
Let $\lambda$ be a polarization of $\et$ and let $H=\lambda_0^{-1}\lambda\in \End_K(\et)$.
Then the polarization $\hat{f^{-1}}\lambda f^{-1}$
of $A$ is defined over $k$ if and only
if $H=\alpha^{\dagger}H\alpha$.
\end{proposition}

\begin{proof}
The polarization will descend to $A$ if and only if it is fixed
by the action of $\sigma$, that is, if and only if
$$
\left(\hat{f^{-1}}\lambda f^{-1}\right)^\sigma
=
\hat{f^{-1}}\lambda f^{-1}
.$$
Multiplying by $f^\sigma$ on the right and by
$\hat{f^\sigma}$ on the left, we find that this
condition is equivalent to
$$
\lambda
=
\hat{f^{-1}f^\sigma}\lambda f^{-1} f^\sigma.$$
This translates into the statement that $H = \alpha^\dag H \alpha$.
\end{proof}

\section{Jacobians in isogeny classes of simple supersingular surfaces}
\label{three}

In this section we will prove Theorem~\ref{T:part2}.  The techniques we use 
depend on whether or not the base field $k=\fq$ has even degree over its prime
field, so we consider these cases in two separate subsections.
Throughout this section we will let $K$ denote an algebraic closure of~$k$,
and we will always assume that the characteristic of $k$ is greater than $3$.

\subsection{The case $q$ a square}
We first show how certain Weil polynomials can be produced by 
Proposition~\ref{P:charpoly}.   We begin with a simple observation:
Suppose $u$ is an automorphism of a hyperelliptic curve $C$, and let $\iota$ be
the hyperelliptic involution of~$C$.  Then $u$ induces an automorphism $u'$ of
the genus-$0$ curve $C/\langle\iota\rangle$, and the order of $u'$ is equal to
the order of $u$ unless $\iota\in\langle u\rangle$, in which case the order of
$u'$ is half that of~$u$.

\begin{proposition}
\label{P:table}
Let $k=\fq$ be a finite field {\rm(}of characteristic at least $5$\/{\rm)} that has
even degree over its prime field.  Let $C$ be a supersingular genus-$2$ 
curve over $k$ such that the Frobenius endomorphism $\pi$ of the 
Jacobian $J$ of $C$ is equal to the integer $\epsilon \sqrt{q}$,
where $\epsilon=\pm 1$.  Let $u$ be a geometric automorphism of $C$ and let
$u'$ be the induced automorphism of $\BP^1$.  Let $n$ and $n'$ be the orders 
of $u$ and $u'$, respectively, and let $C'$ be the twist of $C$ determined
by~$u$.  Then the pair $(n,n')$ appears in the left column of 
Table~\ref{Table:CharPolys}, and the Weil polynomial of $C'$ is 
$q^2 f_{n,n'}(\epsilon x/\sqrt{q})$, where $f_{n,n'}$ is the polynomial 
appearing in the right column.
\end{proposition}

\begin{table}
\begin{center}
\as{1.2}
\begin{tabular}{|c|c|}
\hline
$(n,n')$ & $f_{n,n'}$ \\ \hline\hline
$(1,1)$  & $(x - 1)^4$                    \\ \hline
$(2,1)$  & $(x + 1)^4$                    \\ \hline
$(2,2)$  & $(x - 1)^2(x + 1)^2$           \\ \hline
$(3,3)$  & $(x^2 + x + 1)^2$              \\ \hline
$(4,2)$  & $(x^2 + 1)^2$                  \\ \hline
$(5,5)$  & $x^4 + x^3 + x^2 + x + 1$      \\ \hline
$(6,3)$  & $(x^2 - x + 1)^2$              \\ \hline
$(6,6)$  & $(x^2 - x + 1)(x^2 + x + 1)$   \\ \hline
$(8,4)$  & $x^4 + 1$                      \\ \hline
$(10,5)$ & $x^4 - x^3 + x^2 - x + 1$      \\ \hline
\end{tabular}
\end{center}
\vspace{1ex}
\caption{Characteristic polynomials of certain 
automorphisms of supersingular Jacobians.}
\label{Table:CharPolys}
\end{table}

\begin{proof}
Igusa~\cite[\S8]{igusa} computed the groups that can occur as the reduced 
automorphism groups of hyperelliptic curves.  Looking at Igusa's list, 
we see that $n'$ must be an element of $\{1,2,3,4,5,6\}$.  Since 
$n$ is equal to either $n'$ or $2n'$, we see that the left column includes 
all possibilities except for $n=4$, $n'=4$ and $n=12$, $n'=6$.  These
cases can be excluded by computing the automorphism groups of Igusa's
curves with many automorphisms.

Now Proposition~\ref{P:table} will follow from Proposition~\ref{P:charpoly},
provided that we can show that the characteristic polynomial $g$ of the 
automorphism $u^*$ of $J$ is equal to the polynomial $f$ associated 
to the pair $(n,n')$.  Four facts will be very helpful in our proof that~$g = f$: 
\begin{enumerate}
\item The gcd of $g$ and $x^n-1$ does not divide $x^m-1$ for any $m<n$.
\item If $n=2n'$, the gcd of $g$ and $x^{n'}+1$ does not divide 
      $x^m+1$ for any $m<n'$.
\item If $n=n'$, the gcd of $g$ and $x^n-1$ divides no polynomial
      of the form $x^m + 1$.
\item The constant term of $g$ is $1$.
\end{enumerate}
The first three facts follow from the definitions of $n$ and $n'$ and from the
fact that $u^*$ satisfies the relevant gcd in each case.  The fourth fact 
holds because $u^*$ is an automorphism of the polarized Jacobian, and so 
its product with its Rosati involute is equal to $1$.

These four facts allow us to determine $g$ in all cases, except when $n=3$ 
or $n=6$.  Consider the case when $n=3$.  Then $g$ must be either
$$(x-1)^4, \quad (x-1)^2(x^2+x+1), \text{\quad or\quad} (x^2+x+1)^2.$$
The first possibility is eliminated by fact (1) above. 
If $g$ were equal to the second polynomial, then Proposition~\ref{P:charpoly} 
would show that the Weil polynomial of $C'$ would be
$$(x^2 - 2sx + q)(x^2 + sx + q)$$
for some $s$ with $s^2=q$.  But Theorem~\ref{T:supersingular} 
shows that there are no curves with such a Weil polynomial. 
It follows that $g = (x^2+x+1)^2$, as claimed in the table.
The cases with $n=6$ follow in a similar way.
\end{proof}

\begin{remark}
Suppose $u$ is an order-$n$ automorphism of a genus-$g$ hyperelliptic curve
over an arbitrary algebraically-closed 
field, and let $n'$ be the order of  the automorphism 
of $\BP^1$ induced by~$u$.  In future work, we will show that if $g$
is even or if $n$ is odd, then the characteristic polynomial of $u^*$
is completely determined by the values of $g$, $n$, and~$n'$.
This result provides a much more conceptual way of deriving
Table~\ref{Table:CharPolys}.
\end{remark}

We continue to let $E$ denote an elliptic curve over $\fp$ with trace~$0$.
The following propositions will help us detect the existence of
Jacobians among the supersingular surfaces $A$ that
have respectively $a(A)=2$ or $a(A)=1$.

\begin{proposition}
\label{P:critcurve}
Let $C$ be a supersingular genus-$2$ curve over $k$ whose Jacobian
has $a$-number~$2$.  If $C'$ is a twist of $C$, then the Weil polynomial
of $C'$ is $q^2 f(x/\sqrt{q})$ for some $f = f_{n,n'}$ from 
Table~\ref{Table:CharPolys}.  Furthermore, $C$ has a twist with Weil
polynomial $q^2 f_{n,n'}(x/\sqrt{q})$ if and only if there is a geometric
automorphism $u$ of $C$ of order $n$ that induces an automorphism of 
order $n'$ on the projective line.
\end{proposition}

\begin{proof}
Since the principally polarized surface $\Jac C$ has $a$-number~$2$,
it is geometrically isomorphic to $(E \times E,\lambda)$ for
some principal polarization $\lambda$. Now, $(E \times E,\lambda)$ is
defined over $\mathbb{F}_{p^2}$, hence over $\mathbb{F}_q$, and it is
$k$-isomorphic to the canonically polarized Jacobian of a curve $C_0$
defined over $k$. By Torelli's theorem, $C_0$ is a twist of~$C$. 
Replacing $C_0$ with its quadratic twist, if necessary, we may assume
that the Frobenius on $\Jac C_0$ acts as $\sqrt{q}$.

The proposition is now a direct consequence of Proposition~\ref{P:table}
because $C$ and $C_0$ have the same set of twists, the
$q$-power Frobenius endomorphism $\pi_0$ of the Jacobian of $C_0$
satisfies $\pi_0=\sqrt{q}$, and any isomorphism between the curves
$C$ and $C_0$ induces an isomorphism between their geometric
automorphism groups that identifies the hyperelliptic involutions.
\end{proof}

\begin{proposition}
\label{P:nonprin} Let $\calA$ be the isogeny class $\calA_{(0,2q)}$
{\rm(}respectively $\calA_{(0,0)}$, respectively
$\calA_{(0,-q)}${\rm)}, and let $P \in\Z[x]$ be the polynomial
$(x^2+1)^2$ {\rm(}respectively $x^4+1$, respectively
$x^4-x^2+1${\rm)}. Then there exists an $H\in\Lambda^{\nprc}$ with
an automorphism $\alpha\in\Aut(H)$ such that $P(\alpha)=0$ if and
only if there exists a curve $C$ over $k$ whose Jacobian lies in
the isogeny class $\calA$ and has $a$-number~$1$.
\end{proposition}

\begin{proof}
Suppose $H$ is an element of $\Lambda^\nprc$ and $\alpha$ is an
element of $\Aut(H)$ with $P(\alpha)=0$.  Let $A$ be the $K/k$-twist
of $\et$ determined by $\alpha$. Then by
Proposition~\ref{P:charpoly} we see that $A$ lies in $\calA$. By
Propositions~\ref{P:avalue}, \ref{P:subgroups},
and~\ref{P:goodsubgroups} there is a $k$-rational
$\alpha_p$-subgroup $G$ of $A$ such that $a(A/G)=1$. Now we apply
Proposition~\ref{P:nonprpol} to the degree-$p$ map
$$
\varphi\colon \et\mapright{\sim} A\longrightarrow A/G.
$$
There is a principal polarization $\lambda$ of $A/G$ whose pullback by
$\varphi$ is the degree-$p^2$ polarization $\lambda_0H$ of $\et$,
where $\lambda_0$ is the product principal polarization on~$\et$.
Proposition~\ref{P:polarizations} shows that the pullback of $\lambda$
to $A$ is defined over $k$; hence, $\lambda$ is defined over $k$.
The principally polarized variety $(A/G,\lambda)$ is not geometrically a
product of elliptic curves (because $a(A/G)=1$), so it is the Jacobian of
a curve.

Conversely, suppose that $C$ is a curve over $k$ whose Jacobian
$J$ has $a$-number $1$ and belongs to~$\calA$. Consider the quotient
$J/\alpha_p$ of $J$ by its unique $\alpha_p$-subgroup,
and let $f$ be a geometric isomorphism
$\et\longrightarrow J/\alpha_p$ (which exists by Proposition~\ref{P:avalue}).
Apply Proposition~\ref{P:nonprpol} to the degree-$p$ map
$$\varphi\colon\et\mapright{f} J/\alpha_p\longrightarrow J,
$$
where the rightmost map is the dual isogeny of the canonical projection.
There is some $H\in\Lambda^{\nprc}$ uniquely associated to
the canonical polarization $\theta$ of $J$. Since $\theta$ is defined
over~$k$, Proposition~\ref{P:polarizations} shows that the automorphism
$\alpha=f^{-1}f^{\sigma}$ lies in $\Aut(H)$.
Finally, by Proposition~\ref{P:charpoly} the characteristic polynomial of
$\alpha$ is determined by the class $\calA$ as indicated in the
statement of the proposition.
\end{proof}

Now we proceed to the proof of Theorem~\ref{T:part2} in the case that $q$
is a square.  Consulting Table~\ref{Table:simpleSS}, we see that we must
show that there are Jacobians in $\calA_{(0,0)}$ when $p\not\equiv 1 \bmod 4$,
that there are Jacobians in $\calA_{(0,2q)}$ when $p\equiv 1 \bmod 4$,
and that when $p\not\equiv 1\bmod 3$ there are Jacobians in $\calA_{(0,-q)}$
if and only if $p\equiv 1\bmod 4$.

\subsubsection*{The isogeny class $\calA_{(0,2q)}$
                when $p\equiv 1 \bmod 4$.}

For $p>5$ we deduce from Theorem~\ref{mass2} the existence of a hermitian form
$H\in\Lambda^\nprc$ that admits an automorphism $\alpha$
satisfying $(\alpha^2+1)^2=0$. By Proposition~\ref{P:nonprin}
there is a Jacobian in the class~$\calA_{(0,2q)}$.

For $p=5$ (or more generally for $p\equiv 5 \bmod 8$), we can use
the curve $C$ given by the equation $y^2=x^5-x$. By~\cite[Prop.~1.12]{IKO},
this curve is supersingular and its Jacobian has $a$-number~$2$.
Moreover, the automorphism $u$ given by
$(x,y)\mapsto(-x,\sqrt{-1}\,y)$ satisfies~$u^2=\iota$.
By Proposition~\ref{P:critcurve}, the Jacobian of some twist of $C$ lies
in~$\calA_{(0,2q)}$.\qed

\subsubsection*{The isogeny class $\calA_{(0,0)}$
when $p\not\equiv 1 \bmod 4$.}
First we consider the case $p\equiv 7\bmod 8$.  Let $C$
be the curve $y^2 = x^5 - x$.  By~\cite[Prop.~1.12]{IKO} we know
this curve is supersingular and its Jacobian has $a$-number~$2$.
Moreover, $C$ has a geometric automorphism $u$ satisfying $u^4=\iota$;
for instance, $(x,y)\mapsto(\zeta^2x,\zeta y)$, where $\zeta$ is a primitive
eighth root of unity. By Proposition~\ref{P:critcurve},
the Jacobian of some twist of $C$ lies in~$\calA_{(0,0)}$.

Next we consider the case $p\equiv 3\bmod 8$.  In this case Theorem~\ref{mass2}
shows that there is a hermitian form $H$ in $\Lambda^\nprc$ that has an
automorphism $\alpha$ whose characteristic polynomial is $x^4 + 1$.
By Proposition~\ref{P:nonprin} there is a Jacobian in the class~$\calA_{(0,0)}$.
\qed

\subsubsection*{The isogeny class $\calA_{(0,-q)}$
                when $p\not\equiv 1 \bmod 3$.}

We must show that there is a Jacobian in this isogeny class if and only
if $p\equiv 5\bmod 12$.

To begin with, we note that Proposition~\ref{P:critcurve} shows that 
there is no curve whose Jacobian lies in $\calA_{(0,-q)}$ and has 
$a$-number $2$. 
On the other hand, we see from Proposition~\ref{P:nonprin} that
there will be a curve $C$ over $k$ whose Jacobian
lies in $\calA_{(0,-q)}$ and has $a$-number $1$
if and only if there exists an $H\in\Lambda^{\nprc}$ for which
there is an $\alpha\in\Aut(H)$ satisfying $\alpha^4-\alpha^2+1=0$.
By Theorem~\ref{mass2}, for $p>5$ this
happens if and only if $p\equiv 5\bmod{12}$.

For $p=5$ we note that the matrices $H$ in the
unique equivalence class in $\Lambda^{\nprc}$ have reduced
automorphism group $\pgl2{\GF_5}$ (see Theorem~\ref{mass1}).  This group
contains an element of exact order $6$ that lifts to an
$\alpha\in\Aut(H)$ that must satisfy $\alpha^6=-1$.
In fact, the reduced characteristic
polynomial of $\alpha$ is a power of the minimal polynomial; hence
$\alpha^6=1$ would imply that $\alpha$ has order  $1$, $2$ or $3$ in the
reduced group. \qed

\begin{remark}
The proof that there is no Jacobian with $a$-number $1$ in $\calA_{(0,-q)}$
is also valid for $p=3$.  There is no element of order $6$ in the reduced
automorphism group of the unique equivalence class in $\Lambda^{\nprc}$;
in fact, by Theorem~\ref{mass1} this reduced group is isomorphic to $A_6$.
\end{remark}


\subsection{The case $q$ not a square}
We see from Table~\ref{Table:simpleSS} that to prove Theorem~\ref{T:part2}
in the case where $q$ is not a square, we must show that
there are Jacobians
in $\calA_{(0,0)}$ when $p\not\equiv1\bmod4$,
in $\calA_{(0,q)}$ when $p\not\equiv1\bmod3$,
in $\calA_{(0,-q)}$ when $p\not\equiv1\bmod3$,
and in $\calA_{(0,-2q)}$ for all $p$.

We begin with a remark about twists.  Suppose that $V$ is a variety
over $k=\fq$ all of whose automorphisms are defined over $\ff2$, and let
$\sigma$ denote the Frobenius automorphism of $K$ over~$k$.
Let $\alpha$ be an automorphism of~$V$.
Then there is a $1$-cocycle from $\Gal(K/k)$ to $\Aut_K(V)$ that
sends $\sigma$ to $\alpha$ if and only if $\alpha\alpha^\sigma$
has finite order in $\Aut_K(V)$.  (Note that the latter condition
is equivalent to $\alpha^\sigma \alpha$ having finite order, and that
the orders of $\alpha^\sigma\alpha$ and $\alpha\alpha^\sigma$ are
equal to one another.)  As always, the $K/k$-twists of $V$
correspond to elements of the pointed cohomology set
$H^1(\Gal(K/k), \Aut_K(V))$.  If $V'$ is a twist of $V$ and
$f:V\to V'$ is a $K$-isomorphism, then $V'$ corresponds to the
class of the cocycle that sends $\sigma$ to $f^{-1} f^\sigma$.

The following proposition, similar to Proposition~\ref{P:critcurve}, 
allows us to construct twists
in certain isogeny classes.

\begin{proposition}
\label{P:critcurve2}
Let $C$ be a supersingular genus-$2$ curve over $k$ such that the
Frobenius endomorphism $\pi$ of the Jacobian satisfies 
$\pi^2=\epsilon q$, where $\epsilon=\pm 1$. Let $u$ be a geometric 
automorphism of $C$ such that $u u^{\sigma}$ has order 
$n \in \{1,2,3,4,6\}$  and let $C'$ be the twist of $C$ 
determined by $u$. Then the Weil polynomial
$x^4+a x^3+ b x^2+ a q x +q^2$ of $C'$ is determined by $n$ as 
follows\/{\rm:}

\begin{center}
\begin{tabular}{|c||c|c|c|c|c|}
\hline
$n$     & $1$                & $2$               & $3$              & $4$     & $6$ \\
\hline
$(a,b)$ & $(0,-2\epsilon q)$ & $(0,2\epsilon q)$ & $(0,\epsilon q)$ & $(0,0)$ & $(0,-\epsilon q)$ \\
\hline
\end{tabular}
\end{center}
\end{proposition}

\begin{proof}
Let $J$ be the Jacobian of $C$.  Then the Jacobian $J'$ of $C'$ is 
the twist of $J$ associated to the automorphism $\alpha=u_*$.
Note that $\alpha \alpha^{\sigma}$ has order $n$.

Let $\pi\in\End_K(J)$ and $\pi'\in\End_K(J')$ be the $q$-power Frobenius
endomorphisms of $J$ and $J'$, respectively and let $f\colon J\to J'$ be
a geometric isomorphism such that $\alpha=f^{-1}f^{\sigma}$. 
The condition on $\pi$ implies that $J$ splits over the quadratic 
extension of~$k$, and the condition $(\alpha \alpha^{\sigma})^n=1$ 
implies that $f$ is defined over the extension of $k$ of degree~$24$;
in particular, the isogeny class of $J'$ splits over this extension.
Checking the list of supersingular isogeny classes over odd-degree 
extensions of prime finite fields of odd characteristic that split
over the extension of degree $24$ (see~\cite[Thm.~2.9]{mn} 
and~\cite[Table 1]{mn}), we see that the characteristic polynomial
of $\pi'$ is $x^4+bqx^2+q^2$ for some integer $b\in\{0,\pm 1,\pm 2\}$.

The pullback of $\pi'$ by $f$ is $\alpha\pi$.  Since $\pi^2= \epsilon
q$ and $\alpha^\sigma\pi=\pi\alpha$, this implies that
$$(\alpha\alpha^{\sigma})^2+\epsilon b\,\alpha\alpha^{\sigma}+1=0$$
in $\End_K(J)$.  Comparing this identity with
$(\alpha \alpha^{\sigma})^n=1$, we see that $b=-2 \epsilon, 2
\epsilon, \epsilon, 0$ or $-\epsilon$, according to  $n=1,2,3,4$ or $6$.
\end{proof}

\subsubsection*{The isogeny classes $\calA_{(0,\pm q)}$
when $p\not\equiv 1 \bmod 3$.}

Let $C$ be the curve $y^2 = x^6 + 1$ over~$k$, and let $E$ be the
supersingular elliptic curve $y^2 = x^3 + 1$.  The two obvious maps
from $C$ to $E$ show that the Jacobian $J$ of $C$ is $(2,2)$-isogenous
to $\et$ over~$k$, so the Frobenius $\pi$  of $J$ satisfies $\pi^2 = -q$.

Let $\zeta\in K$ be a primitive sixth root of unity,
and let $u$ be the $K$-automorphism
$$(x,y)\mapsto (\zeta/x, y/x^3)$$
of $C$. One checks easily that
$(u^{\sigma}u)(x,y)=(\zeta^4 x,\zeta^3 y)$;
hence,  $(u^{\sigma}u)^3=\iota$, so $u^\sigma u$ has order~$6$
and the Jacobian of some twist of
$C$ lies in  $\calA_{(0,q)}$ by Proposition~\ref{P:critcurve2}.

If $\zeta$ is a primitive cube root of unity,
the same computation shows that  the automorphism
$u^{\sigma}u$ has order $3$, so the
Jacobian of some twist of $C$ lies in  $\calA_{(0,-q)}$.
\qed

\subsubsection*{The isogeny class $\calA_{(0,-2q)}$.}
For this case, we found it simplest to use
a direct construction involving Kani's result (Theorem~\ref{T:Kani})
combined with Galois descent.

Let $F$ be an elliptic curve over $\ff2$ whose $q^2$-Frobenius is equal
to~$q$, and let $F^{(q)}$ be its Galois conjugate over $\fq$.
The $q^2$-Frobenius acts as the identity on the $2$-torsion points
of~$E$, so all of the $2$-torsion points of $F$ are rational
over~$\GF_{q^2}$.  Label the nonzero $2$-torsion points $P$, $Q$, and $R$,
and let $P^{(q)}$, $Q^{(q)}$, and $R^{(q)}$ be the corresponding points
on~$F^{(q)}$.

We can easily produce four maximal isotropic subgroups of $(F\times F^{(q)})[2]$
that are stable under the action of the Galois group of $\GF_{q^2}$ over
$\GF_q$:
\begin{align*}
  \{(0,0), (P,P^{(q)}), (Q,Q^{(q)}), (R,R^{(q)})\}, \\
  \{(0,0), (P,P^{(q)}), (Q,R^{(q)}), (R,Q^{(q)})\}, \\
  \{(0,0), (P,R^{(q)}), (Q,Q^{(q)}), (R,P^{(q)})\}, \\
  \{(0,0), (P,Q^{(q)}), (Q,P^{(q)}), (R,R^{(q)})\}.
\end{align*}
These are the graphs of certain anti-isometries $F[2]\to F^{(q)}[2]$.
But the number of reducible geometric anti-isometries from $F$ to $F^{(q)}$
is equal to half of the number of geometric isomorphisms from $F$ to $F^{(q)}$,
and since we are in characteristic greater than~$3$, there are at most $6$ such
isomorphisms.  Therefore, at least one of the subgroups $G$ listed above comes
from an irreducible anti-isometry, and so there is a curve $C$ over $\GF_{q^2}$
whose Jacobian is isomorphic to $(F\times F^{(q)}) / G$.  Clearly the
polarized Jacobian of $C$ is isomorphic to its Galois conjugate, so $C$ can
be defined over $\GF_q$.  Furthermore, the isogeny
$\Jac C_{\GF_{q^2}}\to F\times F^{(q)}$ descends to give an isogeny from
$\Jac C$ to the restriction of scalars of $E$.  It follows that $\Jac C$
lies in the isogeny class $\calA_{(0,-2q)}.$\qed

\subsubsection*{The isogeny classes $\calA_{(0,0)}$
when $p\not\equiv 1 \bmod 4$.}
We require two separate arguments for this case, one when $p\equiv7\bmod8$
and one when $p\equiv3\bmod8$.

First suppose that $p\equiv7\bmod8$.
Let $C$ be the curve $y^2 = x^5 - x$. By~\cite[Prop.~1.12]{IKO}
and~\cite[Rem.~1.4]{IKO} we see that the Jacobian $J$ of $C$ is
$k$-isogenous to the product of two supersingular curves. Let
$\zeta\in K$ be a primitive eighth root of unity, and let $u$ be
the $K$-automorphism
$$(x,y)\mapsto (\zeta^2/x, \zeta y/x^3)$$
of $C$.  One checks easily that $(u^{\sigma}u)(x,y)=(-x,\zeta^2 y)$,
so that $u^{\sigma}u$ has order $4$ in $\Aut_K(C)$.
Proposition~\ref{P:critcurve2} then shows that the Jacobian
of some twist of $C$ lies in~$\calA_{(0,0)}$.

Now we turn to the case $p\equiv 3\bmod 8$.  Let $E$ be the
supersingular elliptic curve $y^2 = x^3 + x$,
let $i$ be the geometric automorphism $(x,y)\mapsto(-x,\sqrt{-1}\,y)$
of~$E$, and let $\alpha$ be the automorphism
$$\left[\begin{matrix}
0 & 1\\ i & 0
\end{matrix}\right]$$
of $\et$.  Note that $\alpha^\sigma\alpha$ has order $4$
in $\Aut_K(\et)$, so there is a cocycle from $\Gal(K/k)$ to
$\Aut_K(\et)$ that sends $\sigma$ to $\alpha$.
Let $A$ be the twist of $\et$ corresponding to the cohomology
class in $$H^1(\Gal(K/k), \Aut_K(\et))$$ that contains this cocycle .

Let $\pi$ and $\pi_A$ be the $q$-power Frobenius endomorphisms
of $\et$ and $A$, respectively.  Checking the list of supersingular isogeny
classes over odd-degree extensions of finite prime fields of characteristic
at least~$7$ (see~\cite[Thm.~2.9]{mn} or the Appendix),
we see that the characteristic polynomial of $\pi_A$ is $x^4+bx^2+q^2$,
for some integer $b$. Let $f\colon \et\lra A$ be a geometric isomorphism
such that $\alpha=f^{-1}f^{\sigma}$. The pullback of $\pi_A$ by $f$ is
$\alpha\pi$, so $\alpha\pi$ also has characteristic
polynomial $x^4+bx^2+q^2$.  From the equalities
$$
\pi^2=-q, \quad
\alpha^\sigma\pi=\pi\alpha,\text{\quad and\quad}
 (\alpha\alpha^{\sigma})^2=-1,
$$
we see that we must have $b=0$, so $A$ lies in the isogeny
class~$\calA_{(0,0)}$.

Lemma~\ref{L:numbertheory} below shows that there are positive
integers $r$ and $s$ such that $pr^2 - 2s^2 = 1$.  Let $H$
be the $K$-endomorphism of $\et$ given by
$$H:=\left[\begin{matrix}
pr           & s(1+i)\pi \\
-s\pi(1-i) & pr
\end{matrix}\right]\in\Lambda^{\nprc},$$
and let $\lambda=\lambda_0 H$ be the corresponding degree-$p^2$ polarization
on $\et$, where $\lambda_0$ is the product principal polarization on~$\et$.
One checks easily that
$$H = \alpha^\dag H^\sigma \alpha,$$
where $x\mapsto x^\dag$ is the Rosati involution on $\End_K(\et)$
corresponding to the polarization $\lambda_0$, that is,
the conjugate-transpose involution. Arguing as in the proof of
Proposition~\ref{P:polarizations}, we see that
$\lambda$ descends to a polarization on $A$ defined over $k$.

To complete the proof, we need only find a $k$-rational $\alpha_p$-subgroup
$G$ of $A$ such that $a(A/G) = 1$, for then $\lambda$ will descend to $A/G$
by Proposition~\ref{P:nonprpol}, and the geometrically non-split principally
polarized surface $(A/G,\lambda)$ will be a Jacobian.

By Propositions~\ref{P:avalue} and~\ref{P:subgroups} (which is equally valid
for $q$ nonsquare) we need only find $[i:j]\in\pr1(K)\setminus\pr1(\fpt)$
such that $[i:j] = \widetilde{\alpha}[i^\sigma:j^\sigma]$. Arguing
as in the proof of Proposition~\ref{P:goodsubgroups}, we see that this is
always possible if $q>p^2$. Finally, if $q=p$ not all of the $q+1$ solutions
to  $[i:j] = \widetilde{\alpha}[i^\sigma:j^\sigma]$ can be defined over~$\fpt$;
in fact, these solutions are fixed points of
$\widetilde{\alpha^{\sigma}\alpha}$, and this transformation would be the
identity on $\pr1(\fpt)$. Since $\alpha^{\sigma}\alpha\in\Aut(H)$, this
would imply  $\alpha^{\sigma}\alpha=\pm1$ by Proposition~\ref{det1}, in
contradiction with the condition $(\alpha^{\sigma}\alpha)^2=-1$.
\qed

\begin{lemma}
\label{L:numbertheory}
Let $p$ be a prime that is congruent to $3$ modulo~$8$.  Then there are
positive integers $r$ and $s$ such that $pr^2 - 2s^2 = 1$.
\end{lemma}

\begin{proof}
Let $F = \Q(\sqrt{2p})$, and let $\frakp$ be the prime ideal
of $F$ lying over $p$.  Genus theory shows that the class
number of $F$ is odd, and since $\frakp^2 = (p)$ is principal, we
find that $\frakp$ is principal as well, say $\frakp = (t + s\sqrt{2p})$
for integers $t$ and $s$ that we may take to be positive.  Then we have
$t^2 - 2ps^2 = \pm p$, so $t$ must be a multiple of $p$, say $t = pr$.
We see that then $pr^2 - 2s^2 = \pm 1$.  Considering this equation
modulo~$8$, we find that in fact we must have $pr^2 - 2s^2 = 1$.
\end{proof}



\section{Appendix}

For the sake of completeness we outline a step-by-step procedure
that can be used to check whether a given monic quartic polynomial
$f \in \Z[x]$ is the Weil polynomial of a smooth projective genus-$2$ curve
over a finite field $\GF_q$.  Our main theorem tells when the Weil polynomial
for an abelian surface is the Weil polynomial for a Jacobian, so
mostly what we must do is identify the Weil polynomials of abelian
surfaces.  This has been done in other papers (\cite{ruc} and \cite{mn}
for example); we are simply restating these results in a convenient form.

Write $q = p^m$ for a prime $p$.

\begin{itemize}
\item\emph{Step 1\/{\rm:} Check whether $f$ has the right shape for a Weil polynomial.}
The complex roots of the Weil polynomial of an abelian variety over $\GF_q$
have magnitude $\sqrt{q}$ and come in complex conjugate pairs.
A monic quartic polynomial in $\Z[x]$ has this property if and only if
it has the shape
$$f=x^4+a x^3+ b x^2 + q a x+q^2,$$ with
$$|a| \leq 4 \sqrt{q} \text{\quad and\quad}
2|a| \sqrt{q}-2q \leq b \leq \frac{a^2}{4}+2q.$$

\item\emph{Step 2\/{\rm:} Check whether $f$ is the Weil polynomial of an abelian surface.}
Suppose $f$ meets the condition of Step 1, and let
$$\Delta=a^2-4(b-2q)\text{\quad and\quad} \delta=(b+2q)^2-4 q a^2.$$

\begin{itemize}
\item \emph{Ordinary case}\/:  $v_p(b)=0$.

In this case the polynomial $f$ is the Weil polynomial of an ordinary abelian
surface over $\GF_q$. The surface is split or simple according to $\Delta$
being a square in $\Z$ or not.

\item \emph{Mixed case}\/: $v_p(a)=0$ and $v_p(b)>0$.

The polynomial $f$ is the Weil polynomial of a mixed abelian surface if
and only if
$$v_p(b) \geq m/2
\text{\quad and $\delta$ is either $0$ or a non-square in $\Z_p$.}$$
The surface is split or simple according to $\Delta$ being a square in $\Z$
or not.

\item \emph{Supersingular case}\/: $v_p(a)>0$ and $v_p(b)>0$.

The polynomial $f$ is the Weil polynomial of a supersingular split abelian
surface if and only if
$$v_p(a) \geq m/2, \quad v_p(b) \geq m, \text{\quad and\quad} \Delta\mbox{ is a square in }\Z,$$
and moreover, if $q$ is a square and we write $a=\sqrt{q}a'$ and $b=qb'$,
the following two conditions hold:
\begin{align*}
p\not\equiv1\bmod4,& \text{$\quad$ if } b'=2,\\
p\not\equiv1\bmod3,& \text{$\quad$ if } a'\not\equiv b'\bmod2.
\end{align*}

The polynomial $f$ is the Weil polynomial of a simple supersingular abelian
surface if and only if $(a,b)$ belongs to the following list : \\

\begin{center}
\as{1.1}
\begin{tabular}{|c|l|}
\hline
$(a,b)$              & Conditions on $p$ and $q$                       \\
\hline\hline
$(0,0)$              & $q$ is a square and $p\not\equiv1\bmod8$, or    \\
                     & $q$ is not a square and $p \ne 2$               \\
\hline
$(0,-q)$             & $q$ is a square and $p\not\equiv1\bmod{12}$, or \\
                     & $q$ is not a square and $p \ne 3$               \\
\hline
$(0,q)$              & $q$ is not a square                             \\
\hline
$(0,-2q)$            & $q$ is not a square                             \\
\hline
$(0,2q)$             & $q$ is a square and $p \equiv 1 \bmod4$         \\
\hline
$(\pm \sqrt{q},q)$   & $q$ is a square and $ p \not \equiv 1 \bmod5$   \\
\hline
$(\pm \sqrt{2q},q)$  & $q$ is not a square and $p=2$                   \\
\hline
$(\pm 2\sqrt{q},3q)$ & $q$ is a square and $p \equiv 1 \bmod3$         \\
\hline
$(\pm \sqrt{5q},3q)$ & $q$ is not a square and $p=5$                   \\
\hline
\end{tabular}
\vskip 1em
\end{center}

\end{itemize}

\item \emph{Step 3\/{\rm:} Apply Theorem~\ref{T:main}}.
If $f$ is the Weil polynomial of an abelian surface over $\GF_q$, one applies
Theorem~\ref{T:main} to determine if it is the Weil polynomial of a genus-$2$
curve.  Note that in the split case, $\Delta$ is a square in $\Z$ and the
polynomial $x^2+a x+(b-2q)$ has two roots $s,t \in \Z$, which are the traces
of Frobenius of the corresponding elliptic curves.

\end{itemize}


\end{document}